\documentclass[aap,MSNbibl,citesort,dvips]{arximspdf}

\doi{10.1214/09-AAP635}
\volume{20}
\issue{2}
\pubyear{2010}
\firstpage{722}
\lastpage{752}

\makeatletter
\DeclareMathAlphabet\mathcaligr{OMS}{cmsy}{m}{n}
\newcommand{\eqref}[1]{(\ref{#1})}
\def\1{\mathbf{1}}

\newproclaim{rem}{Remark}[section]
\newtheorem{prop}{Proposition}[section]
\newproclaim{defi}{Definition}[section]
\newtheorem{Th}{Theorem}[section]
\newtheorem{coro}{Corollary}[section]
\newtheorem{lem}{Lemma}[section]

\makeatother

\begin{document}
\begin{frontmatter}

\title{The limiting move-to-front search-cost in law of~large~numbers
asymptotic regimes}
\runtitle{L.L.N. for random partitions}

\begin{aug}
\author[A]{\fnms{Javiera} \snm{Barrera}\thanksref{t1,t3}\ead[label=e1]{javiera.barrera@uai.cl}}
and
\author[B]{\fnms{Joaqu\'{\i}n} \snm{Fontbona}\thanksref{t2,t3}\ead[label=e2]{fontbona@dim.uchile.cl}\corref{}}
\pdfauthor{Javiera Barrera, Joaquin Fontbona}
\runauthor{J. Barrera and J. Fontbona}
\affiliation{Universidad Adolfo Ib\'{a}\~{n}ez and Universidad de Chile}
\address[A]{
Escuela de Ingenier\'{\i}a y Ciencias\\
Universidad Adolfo Ib\'{a}\~{n}ez\\
Avda. Diagonal Las Torres 2640\\
Edificio C\\
Pe\~{n}alolen\\
Chile\\
\printead{e1}} 
\address[B]{
Departamento de Ingenier\'{\i}a Matem\'{a}tica and\\
\quad Centro de Modelamiento Matem\'{a}tico\\
DIM-CMM, UMI(2807) UCHILE-CNRS\\
Universidad de Chile\\
Casilla 170-3, Correo 3\\
Santiago\\
Chile\\
\printead{e2}}
\end{aug}
\thankstext{t1}{Supported by Fondecyt Project 1060485.}
\thankstext{t2}{Supported by Fondecyt Project 1070743.}
\thankstext{t3}{Supported by BASAL-Conicyt CMM-U de Chile and
Millennium Nucleus
Information and Randomness ICM P04-069.}
\received{\smonth{11} \syear{2006}}
\revised{\smonth{8} \syear{2009}}

%
\begin{abstract}
We explicitly compute the limiting transient distribution of the
search-cost in the move-to-front Markov chain
when the number of objects tends to infinity, for general families of
deterministic or random request
rates. Our techniques are based on a ``law of large numbers for
random partitions,'' a~scaling limit that allows us to exactly
compute limiting expectation of empirical functionals of the
request probabilities of objects. In particular, we show that the
limiting search-cost can be split at an explicit deterministic
threshold into one random variable in equilibrium, and a second
one related to the initial ordering of the list. Our results
ensure the stability of the limiting search-cost under general
perturbations of the request probabilities. We provide the
description of the
limiting transient behavior in
several examples where only the stationary regime is known, and
discuss the range of validity of our scaling limit.
\end{abstract}

%
\begin{keyword}[class=AMS]
\kwd[Primary ]{60B10}
\kwd{68W40}
\kwd[; secondary ]{68P10}.
\end{keyword}
\begin{keyword}
\kwd{Move-to-front rule}
\kwd{search-cost}
\kwd{law of large numbers}
\kwd{propagation of chaos}.
\end{keyword}

\end{frontmatter}

\section{Introduction}

We consider the
search-cost process in the move-to-front (MtF) Markov chain.
A finite set of objects labeled $1,\ldots, n$ is dynamically
maintained as a serial list, and objects are requested at random
instants with a given probability $p^{(n)}_i $, $i=1,\ldots,n$.
Instantaneously after request, an object is moved to the front
of the list, while the relative order of the other objects is left
unchanged. The search-cost at a given instant is defined as the
position in the list of the next requested object.

The exact and limiting behaviors of the move-to-front rule have
received much attention in the computer science and discrete
probability literature since the 1960s (see Fill \cite{FillMtF}
and Jelenkovi\'{c} \cite{Jelen} for historical references).
For a fixed and finite number of objects, the search-cost distribution
has been
studied by Fill \cite{FillMtF96,FillMtF}, Fill and Holst
\cite{FillHolst96} and Flajolet, Gardy and Thimonier
\cite{FlajoletGardyThimonier} for different deterministic request
probability vectors $\mathbf{p}^{(n)}= (p^{(n)}_1, \ldots,
p^{(n)}_n )$. To our knowledge, the limiting search-cost
distribution as the number of objects goes to infinity is known
only for the stationary regime of the MtF Markov chain. This
problem was first studied by Fill \cite{FillMtF96} for several
types of deterministic request probabilities $\mathbf{p}^{(n)}$,
and later by Barrera, Huillet and Paroissin
\cite{BarreraHuilletParoissinMtF,BarreraHuilletParoissinMtF2}
and Barrera and Paroissin \cite{BarreraParoissinMtF} for random
request probabilities $\mathbf{p}^{(n)}$ defined by normalized
samples of positive i.i.d. random variables. A different approach
was adopted by Jelenkovi\'{c} \cite{Jelen}, who considered various
scaling limits for the stationary search-cost when the optimal
static distribution of objects in the list is specified.

In this article we explicitly compute the
limiting law of the transient search-cost as the number $n$ of
objects tends to infinity, for a large class of deterministic or
random request probabilities. This class includes several
previously considered cases. To be more precise, let
$w_1^{(n)},\ldots, w_n^{(n)}\geq0$ be deterministic or random real
numbers and consider the probability vector
$\mathbf{p}^{(n)}= (p_i^{(n)} )_{i=1}^n$ defined by
%
\begin{equation}\label{randompartition}
p_i^{(n)}=\frac{w_i^{(n)}}{\sum_{j=1}^n w_j^{(n)}}.
\end{equation}
We call $p_i^{(n)}$ the ``popularity'' of object $i$.
Let us
denote by ${\mathcaligr P}(\mathbb{R}_+)$ the space of Borel probability
measures in $\mathbb{R}_+$. We consider request probabilities
$\mathbf{p}^{(n)}$ of this type that
exhibit a weak law of large numbers behavior.
That is, we assume that as~$n$ goes to infinity, the empirical
measure
\[
\hat{\nu}^{(n)}:=\frac{1}{n}\sum_{i=1}^n\delta_{w_i^{(n)}}\in
{\mathcaligr P}(\mathbb{R}_+),\qquad   n\in\mathbb{N},
\]
converges in distribution [as a random element of the Polish space
${\mathcaligr P}(\mathbb{R}_+)$] to a deterministic limit $P\in
{\mathcaligr
P}(\mathbb{R}_+)$. Moreover, we assume that $\int x \hat{\nu}^{(n)}(dx)$
converges to $ \mu:= \int x
P(dx)\in(0,\infty)$ in distribution.
(These conditions hold if, for instance, $w_i^{(n)}=w_i$,
$i=1,2,\ldots,$ are i.i.d. random variables with finite mean $\mu$.)

Denote by $\nu^{(n)}$ the empirical measure of the scaled vector
$(n\mu) \mathbf{p}^{(n)}$. As we shall see, the computation of the
search-cost distribution will involve continuous functionals of
$\nu^{(n)}$. On the other hand, we will show that the sequence $\nu^{(n)}$ shares
the same law of large numbers behavior and
limit of $\hat{\nu}^{(n)}$. This fact is what we call ``law of
large numbers for random partitions of the interval.'' With these
elements, we will be able to explicitly compute for all $t>0$ the
limiting distribution of the (suitably normalized) transient
search-cost, in terms of the limiting probability measure $P$.

We will
show the existence of an explicit deterministic threshold,
depending only on $P$ and $t$, such that the limiting transient
and stationary search-costs have the same distribution in the
event they fall below it. The limiting transient search-cost
restricted to that event will be called ``equilibrium
part'' of the transient search-cost, and its distribution (as well as the
limiting stationary one) will depend only on the law $P$.
Alternatively, we will call ``out-of-equilibrium part'' of the
transient search-cost its restriction to the complementary event.
The asymptotic behavior of the out-of-equilibrium part will
depend on $P$, but also on the relative order of popularities in
the list at $t=0$. Concretely, three situations will be studied:
(1)~no information at all about the request probabilities is
available at the beginning; (2)~objects in the list are known to be
initially arranged in decreasing order\break of popularity; and (3)~objects
are known to be initially arranged in increasing order of
popularity. In all three cases, we find an explicit expression
for the limiting law~of the out-of-equilibrium part in terms of
$P$.

As a consequence, we obtain an upper bound for the total
variation distance between the transient and the stationary
limiting search-cost distributions. The distance from equilibrium
at time $t$ turns out to behave like $O(\int x e^{-tx}P(dx))$ in
each of the three situations we consider. The techniques we
introduce also show that the limiting transient search-cost
distribution is stable under perturbations of the request
probabilities preserving both the limiting law $P$ and the initial
relative order in the list.

The rest of this article is organized as follows.
In Section~\ref{sec2} we define the MtF Markov chain in continuous time and
the associated search-cost process following the lines of Fill and
Holst \cite{FillHolst96}. We also recall some nonasymptotic
results about its law. In Section~\ref{sec3} we state our main
results, namely Theorems~\ref{theo:principaleq} and~\ref{theo:principalneq}, which provide the limiting expressions
for the Laplace transforms of the two components of the suitably
normalized transient search-cost. We deduce from them the limiting
search-cost law in terms of $P$ and its Laplace transform. We also
present examples, we discuss connections with the
Persistent-Caching-Algorithm introduced by Jelenkovi\'{c} and
Radovanovi\'{c} \cite{JR} and we compare our type of law of large
numbers asymptotic with the fluid limit considered by Jelenkovi\'{c}
\cite{Jelen}. Furthermore, we prove stochastic order relations
between the search-costs in the three situations (regarding the
initial ordering) that are considered. In Section~\ref{sec4}, we
prove our law of large numbers for random partitions, we discuss
its connection with the propagation of chaos property arising in
the probabilistic study of mean field models, and we use those
ideas to prove Theorems~\ref{theo:principaleq} and~\ref{theo:principalneq}. In the last section, we discuss the scope
of application and the limitations of our techniques.

Let us establish some notation.
In the sequel, $\to^d$ means convergence in
distribution, and the notation $\Longrightarrow$ stands for weak
convergence of probability measures. We denote by $\delta_x$ the
Dirac mass at some point $x$. The convention $\frac{0}{0}:=0$ is
adopted throughout.

\section{Preliminaries and notation}\label{sec2}

Consider a list of $n$ objects labeled $\{1,\ldots,n\}$ and a
permutation $\pi$ of
$\{1,\ldots,n\}$. Assume that at time $t=0$, object $i$ is at
position $\pi(i)$ of the list. Then, objects are requested at
random instants $t>0$ after which the list is instantaneously
modified,
by placing the requested object on its top. This is
the MtF rule. It is customary to assume that
different objects are requested at random instants given by
independent standard Poisson processes in the line. Let~$w_i$
denote the intensity at which object $i$ is requested. The total
number\vspace*{1pt} of requests up to time $t$ defines a Poisson process, say
$\tilde{N}_t$, of rate
\[
w=\sum_{j=1}^n w_j.
\]
By the strong
Markov property, the probability that object $i$ is requested at a
given arrival time of $\tilde{N}_t$ is
\[
p_i:=\frac{w_i}{w}.
\]
We call this quantity the ``popularity'' of object $i$.

We shall in the sequel work with the time-changed process
\[
N_t:=\tilde{N}_{t/w}
\]
and its requests instants. This is
the time scale considered, for instance, \mbox{in \cite{FillHolst96,Bodell},} and in our case it will simplify the asymptotic
analysis by
keeping a constant (unitary) total rate of requests.
The request rate $w_i$ of object $i$ becomes $p_i$ in the new time
scale, but its popularity remains unchanged.

\begin{rem} Nevertheless, our statements will rely on
hypotheses made on the parameters $(w_i)$ and will be interpreted
also in the original time scale (see Remark~\ref{originaltime}
below).
\end{rem}

We denote by
\[
S^{(n,i)}(t)
\]
the position in the list at time $t$ of object
$i$, and by $I_k\in\{1,\ldots,n\}$ the $k$th requested object
(in chronological order). Thus, the label of the first object
requested in the time interval $[t,+\infty)$ is $I_{N_{t^-}+1}$.
We are interested in the search-cost of that object. That is, in the random
variable defined by
\[
S^{(n)}(t):=\sum_{i=1}^n S^{(n,i)}(t)\1_{\{I_{N_{t^-}+1}=i\}}.
\]
Notice that although the processes $S^{(n,i)}(t)$ are left
continuous, $S^{(n)}(t)$ is right continuous since the list is
modified instantaneously after each request.

We will further need the following notation:
\begin{longlist}[$\bullet$]
\item[$\bullet$] $R_t$ is the subset of $\{1,\ldots,n\}$ consisting of objects that
have been
requested at least once in the time interval $[0,t[$.
\item[$\bullet$] We decompose the search-cost $S^{(n)}(t)$ into two random
variables:
\[
S^{(n)}(t)=S^{(n)}_{\mathbf{e}}(t)+S^{(n)}_{\mathbf{o}}(t),
\]
where
\[
S_{\mathbf{e}}^{(n)}(t):=S^{(n)}(t)\1_{\{I_{N_{t^-}+1} \in R_t\}}
\]
and
\[
S_{\mathbf{o}}^{(n)}(t):=S^{(n)}(t)\1_{\{I_{N_{t^-}+1} \notin
R_t\}}.
\]
\end{longlist}

Thus, $S_{\mathbf{e}}^{(n)}(t)$ is the search-cost of the
requested object if it has been requested at least once in
$[0,t[$, and it is $0$ otherwise. $S_{\mathbf{o}}^{(n)}(t)$ is
defined conversely. The subscripts ${\mathbf{e}}$ and
${\mathbf{o}}$ respectively stand for ``equilibrium'' and ``out of
equilibrium.'' This decomposition and notation are inspired in
Fill's work \cite{FillMtF}, where a coupling was introduced which
simultaneously updates the list in stationary regime and an
arbitrary second list. In that coupling, each object had the
same search-cost in the two lists after its first request.

\begin{rem}\label{SC} Notice that an object has been requested before
time $t$ if
and only if it stands at one of the first $|R_t|$ positions in
the list. Therefore, we have
\[
\bigl\{S_{\mathbf{e}}^{(n)}(
t)>0\bigr\}=\bigl\{S^{(n)}(t)\leq|R_t|\bigr\}.
\]
\end{rem}

The next result will be used in the sequel.

\begin{prop}\label{prop:bodell}
Let $\pi(i)$ be the position in the list of item $i$ at time $0$.
For given real parameters $q_1,\ldots,q_n\in[0,1]$ let
$(B_1(q_1)\ldots,B_n(q_n))$ denote a vector of $n$ independent
Bernoulli random variables of such parameters.
\begin{longlist}
\item[(a)] For all
$k,i\in\{1,\ldots,n\}$,
\[
\mathbb{P}\bigl\{S^{(n,i)}(t)=k,i\in R_t\bigr\}= \int_0^t p_i e^{-p_i u}
\mathbb{P}
\bigl\{J_{\mathbf{e}}^{(n)}(u)=k\bigr\}\, du ,
\]
where $ {J_{\mathbf{e}}^{(n)}(u)=^d \sum_{j=1, j\not=
i}^n B_j(1-e^{-p_j u}) }$.
\item[(b)] For all
$k,i\in\{1,\ldots,n\}$,
\[
\mathbb{P}\bigl\{S^{(n,i)}(t)=k,i\notin R_t\bigr\}= \mathbb{P}
\bigl\{J_{\mathbf{o}}^{(n)}(t)=k\bigr\}e^{-p_i t},
\]
where $ {J_{\mathbf{o}}^{(n)}(t)=^d \sum_{j=1,j\not=
i}^n B_j (1-e^{-p_j t}\mathbf{1}_{{\pi(i)<\pi(j)}} ) }$
and $\mathbf{1}_{{\pi(i)<\pi(j)}}=1$ if object $i$ precedes $j$
in the initial permutation or $0$ otherwise.
\item[(c)] For all
$k\in\{1,\ldots,n\}$,
\[
\mathbb{P}\bigl\{S_{\mathbf{e}}^{(n)}(t)=k\bigr\}=\sum_{i=1}^n \int_0^t p_i^2
e^{-p_i u} \mathbb{P}\bigl\{J_{\mathbf{e}}^{(n)}(u)=k\bigr\}\, du.
\]
\item[(d)]For all
$k\in\{1,\ldots,n\}$,
\[
\mathbb{P}\bigl\{S_{\mathbf{o}}^{(n)}(t)=k\bigr\}=\sum_{i=1}^n p_i \mathbb{P}
\bigl\{J_{\mathbf{o}}^{(n)}(t)=k\bigr\}e^{-p_i t}.
\]
\end{longlist}
\end{prop}

\begin{pf} The proof of relations {(a)} and {(b)} can be
deduced from Proposition~2.1 in \cite{FillHolst96}.
The basic ideas are to condition in the last instant $u\in\,\,]0,t]$
where object~$i$ has been requested and to consider the Poisson
point process in reversed time starting from $t$
(see Theorem 2.3.1.3. and Corollary 2.3.1.6
in \cite{Bodell} for a complete proof).

Relation {(c)} [resp. {(d)}] follows easily from {(a)}
[resp. {(b)}], thanks to independence of the events
$\{S^{(n,i)}(t)=k,i\in R_t\}$ [resp. $\{S^{(n,i)}(t)=k,i\notin
R_t\}$] and \mbox{$\{I_{N_{t^-}+1}=i\}$.}
\end{pf}

\section{Main statements, examples and consequences}\label{sec3}

For each $n$ we next consider a random or deterministic
vector of nonnegative real numbers
\[
\mathbf{w}^{(n)}= \bigl(w_1^{(n)}, \ldots, w_n^{(n)} \bigr).
\]
Then, \textit{conditionally} on $\mathbf{w}^{(n)}$, we define the
MtF Markov chain and its search-cost $S^{(n)}(t)$ in the
same way as was done in the previous section for deterministic
request rates. Recall that the process $S^{(n)}(t)$ refers to the
time-scale at which requests arrive at rate $1$.

\begin{rem} By Proposition \ref{prop:bodell}, the law of $S^{(n)}(t)$
conditional on $\mathbf{w}^{(n)}$ depends on that vector only
through the popularities
\[
\mathbf{p}^{(n)}= \bigl(p_1^{(n)},\ldots,p_n^{(n)} \bigr)
\]
defined as in \eqref{randompartition}.
\end{rem}

Let us recall the result
obtained in \cite{BarreraHuilletParoissinMtF2} for request
probability vectors given by normalized samples of positive i.i.d.
random variables.

\begin{Th}\label{javiera}
Let $(w_i)_{i\in\mathbb{N}}$ be an i.i.d. sequence of nonnegative
random variables with finite mean $\mu$ and Laplace transform
$\phi(t)$ and, for each $n\in\mathbb{N}$, take
$\mathbf{w}^{(n)}=(w_i)_{i=1}^n$.

Let $S^{(n)}(\infty)$ be a random variable defined,
conditionally on $\mathbf{w}^{(n)}$, as the search-cost associated
with the MtF Markov chain in stationary regime. Then, when $n\to
\infty$, we have the convergence
\[
\frac{S^{(n)}(\infty)}{n}\rightarrow^d S^{(\infty)},
\]
where $S^{(\infty)}$ is a random variable in $[0,1]$ with density
given by
\[
f_{S^{(\infty)}}(x)=-\frac{1}{\mu}\frac{\phi'' (\phi
^{-1}(1-x) )}{\phi' (\phi^{-1}(1-x) )}
\1_{[0,1-\mathbf{p}_0]},
\]
and $\mathbf{p}_0=\mathbb{P}(w_i=0)$.
\end{Th}

The proof of Theorem~\ref{javiera} relied on Laplace integral
techniques. Our goal now is to describe the behavior of a suitable
normalization of the random variable $S^{(n)}(t)$ when $n$ goes to
$\infty$. Furthermore, we will do this under an assumption
naturally generalizing that of
\cite{BarreraHuilletParoissinMtF2}.

\begin{defi}[(Condition LLN-$P$)]
We say that a sequence of (random or deterministic) vectors
$\mathbf{w}^{(n)}= (w_1^{(n)}, \ldots, w_n^{(n)} )_{n \in
\mathbb{N}}$ satisfies a law of large numbers with limiting law $P$
(LLN-$P$ for short), if there exist a probability measure $P\in
{\mathcaligr P}(\mathbb{R}_+)$ with finite first moment $\mu\not=0$ and
positive random variables $Z_n$, such that the empirical measures
\[
\hat{\nu}^{(n)}:=\frac{1}{n}\sum_{i=1}^n \delta_{Z_n w_i^{(n)}}
\]
converge in law to $P$, and their empirical means
\[
\frac{1}{n}\sum_{i=1}^n Z_n w^{(n)}_i
\]
converge in law to
$\mu$.
\end{defi}

\begin{rem}
\begin{longlist}
\item[(i)] Condition LLN-$P$ may hold true
with random variables
$Z_n$'s that are not identically equal to $1$ and fail to hold if
one takes $Z_n\equiv1$ [see, e.g., {(c)} below].
\item[(ii)] LLN-$P$ is equivalent to say that the sequence
$(\hat{\nu}^{(n)})$ converges in distribution to the deterministic
value $P$, when seen as random variables in the Polish space $
{\mathcaligr P}_1(\mathbb{R}_+) $ of Borel probability measures with finite
first moment, endowed with the Wasserstein distance $W_1$ (see
Theorem~\ref{teoLLN}).
\item[(iii)] In the time scale originally introduced in Section~\ref{sec2},
the constant $\mu$ may be thought of as the (asymptotic in $n$)
average request rate per object.

\end{longlist}
\end{rem}

Provided that the empirical means converge in law, LLN-$P$
holds in several situations. The following are some examples:

\begin{longlist}
\item[(a)] $w^{(n)}_i=w_i$
for all $n\in\mathbb{N}$, with $(w_i)_{i\in\mathbb{N}}$ an ergodic process
with invariant measure $P$, and $Z_n=1$.
\item[(b)] $(w_1^{(n)},\ldots,w_n^{(n)})$, $n\in\mathbb{N}$,
is an exchangeable and $P$-chaotic vector and $Z_n=1$.
Recall that a random vector $v=(v_i)_{i=1}^n$ in $\mathbb{R}^n$ is said
to be \textit{exchangeable} if the law of $(v_{\sigma(i)})_{i=1}^n$
is the same for any $n$-permutation $\sigma$. The notion of
``$P$-chaotic vector'' is recalled in Section~\ref{sec4}.
\item[(c)] $w_i^{(n)}=i^{\alpha}$ for all $n\in\mathbb{N}$, $\alpha
\in\mathbb{R}
$ and $Z_n=n^{-\alpha}$. Indeed, for any
$\varphi\dvtx\mathbb{R}\to\mathbb{R}$ continuous and bounded, one has
$\frac{1}{n}\sum_{i=1}^n
\varphi (\frac{i^{\alpha}}{n^{\alpha}} )\to\int_0^1
\varphi(x^{\alpha})\,dx$, using the continuity of $x\mapsto
x^{\alpha}$ in $(0,1]$ for any $\alpha$. By the obvious change of
variable we get
%
\begin{equation}\label{Pexamplec}
P(dx)=
\cases{
\dfrac{1}{\alpha}x^{1/\alpha-1}\1_{[0,1]}(x)\,dx,
&\quad if $\alpha> 0$,\vspace*{2pt} \cr
\delta_1(dx),
&\quad if $\alpha= 0$,\vspace*{2pt} \cr
\dfrac{1}{|\alpha|}x^{1/\alpha-1}\1_{[1,\infty)}(x)\,dx,
&\quad if $\alpha<0$.
}
\end{equation}
Thus, one can check that LLN-$P$ holds for these
$w_i^{(n)}$ if and only if $\alpha>-1$.
\item[(d)] Let $q$ be a continuous probability density with compact
support in $[0,c]$, and for each $n\in\mathbb{N}$
define
\[
w^{(n)}_i:=Q(ci/n)- Q(c(i-1)/n),\qquad
i=1,\ldots, n,
\]
where $Q(x)=\int_0^x q(y)\,dy$ is the primitive of $q$. Then, setting
$Z_n:=n/c$, for certain $x_i^{(n)}\in((i-1)/n,i/n]$,
$i=1,\ldots, n$, we have
\[
\frac{1}{n}\sum_{i=1}^n \varphi \bigl(Z_n w_i^{(n)}  \bigr) =
\frac{1}{n}\sum_{i=1}^n \varphi\bigl( q\bigl(c x_i^{(n)}\bigr)\bigr).
\]
Hence, for any
continuous and bounded function $\varphi\dvtx\mathbb{R}\to\mathbb{R}$,
we get
\[
\frac{1}{n}\sum_{i=1}^n \varphi\bigl(Z_n w_i^{(n)}\bigr)\to\frac{1}{c}
\int_0^c \varphi(q( x))\,dx,
\]
and so LLN-$P$ holds in this case with $P(dx)=(\frac{1}{c}\1_{x\in[0,c]}\,dx)\circ q^{-1}$, the push-forward of the
normalized Lebesgue's measure by $q$. (The convergence of empirical
means is in this case trivial.)
\end{longlist}

Case {(c)} above was considered by Fill \cite{FillMtF96}. Case
{(d)} is a particular instance of the ``light tail'' condition
studied in Jelenkovi\'{c} \cite{Jelen}. See the examples below for
more details.

Before stating our main results, notice that from Proposition
\ref{prop:bodell} the law of $S^{(n)}_{\mathbf{e}}(t)$ does not
depend on the initial permutation $\pi$ of the list, whereas that
of $S^{(n)}_{\mathbf{o}}(t)$ does. This is simply due to the fact
that the cost of searching an already requested object does not
depend any more on its initial position. In turn, the value of the
initial permutation $\pi$ remains present in the law of
$S^{(n)}_{\mathbf{o}}(t)$. By this reason, we need to study
separately the two components of the transient search-cost. Let us
state our main theoretical result on the equilibrium part.

\begin{Th}\label{theo:principaleq}
For $\lambda\geq0$ and $t\geq0$, define
\begin{eqnarray*}
A_n(t,\lambda)
&:=&
 \mathbb{E} \bigl(\exp\bigl\{-\lambda
S_{\mathbf{e}}^{(n)}(t)\bigr\}\1_{\{S^{(n)}_{\mathbf{e}}(t)>0
\}} \bigr) \\
&&\hspace*{-10.8pt}= \mathbb{E} \bigl(\exp\bigl\{-\lambda S^{(n)}(t)\bigr\}\1_{\{S^{(n)}(t)\leq
|R_t|\}}  \bigr).
\end{eqnarray*}
Then, if LLN-$P$ holds, we have
\[
\lim_{n\to\infty}A_n(n\mu t,\lambda/ n)= \frac{1}{\mu} \int
_0^t \int_{\mathbb{R}^+} x^2e^{-x u} P(dx) \exp\bigl\{-\lambda\bigl(1-\phi
(u)\bigr) \bigr\}\, du,
\]
where $\phi\dvtx\mathbb{R}_+\to\mathbb{R}_+$ is the Laplace transform
of $P$.
\end{Th}

\begin{rem}\label{recovstat}
The Laplace transform of the limiting stationary search-cost
obtained
in \cite{BarreraHuilletParoissinMtF2} corresponds to the limit
of the latter expression when $t\to\infty$.
\end{rem}

By the above exposed reasons, some asymptotic assumptions on the
initial ordering $\pi$ of the list will be needed in order to
observe a coherent limiting behavior of the out-of-equilibrium
part of the transient search-cost $S^{(n)}_{\mathbf{o}}(t)$.
Notice that any relevant property of $\pi$ can be restated in
terms of the vector of popularities
$\mathbf{p}^{(n)}$, and one can therefore assume
without loss of generality that $\pi$ is equal to the identity
permutation $Id$. We shall explicitly analyze three particular
assumptions on
$\mathbf{w}^{(n)}$ or (equivalently) on $\mathbf{p}^{(n)}$:
\begin{eqnarray*}
&&\mbox{LLN-$P${-ex}: \quad LLN-$P$ holds, $\pi=Id$ and
$\mathbf{w}^{(n)}$ is
exchangeable for each $n\in\mathbb{N}$.}\\
 &&\mbox{LLN-$P^-$: \quad\hspace*{7pt}LLN-$P$ holds, $\pi=Id$ and
$\mathbf{w}^{(n)}$ is decreasing a.s. for each $n\in\mathbb{N}$.}\\
 &&\mbox{LLN-$P^+$: \quad\hspace*{7pt}LLN-$P$ holds, $\pi=Id$ and
$\mathbf{w}^{(n)}$ is increasing a.s. for each $n\in\mathbb{N}$.}
\end{eqnarray*}

Clearly,
the assumption $\pi=Id$ is
superfluous under LLN-$P${-ex}, but we shall adopt
it for notational convenience.
The asymptotic behavior of the out-of-equilibrium part of the
transient search-cost is stated in the following:

\begin{Th}\label{theo:principalneq}
For $\lambda\geq0$ and $t\geq0$, define
\begin{eqnarray*}
B_n(t,\lambda)
&:= & \mathbb{E} \bigl(\exp\bigl\{-\lambda
S_{\mathbf{o}}^{(n)}(t)\bigr\}\1_{\{S^{(n)}_{\mathbf{o}}(t)>0\}} \bigr)\\
&&\hspace*{-10.8pt}=
\mathbb{E} \bigl(\exp\bigl\{-\lambda S^{(n)}(t)\bigr\}\1_{\{
S^{(n)}(t)>|R_t|\}
} \bigr)
\end{eqnarray*}
and define $\phi\dvtx\mathbb{R}_+\to\mathbb{R}_+$ as before. Then,
$L(\mu,t,\lambda):=\lim_{n\to\infty} B_n(n\mu t,\lambda/n)$
exists in the following cases:
\begin{longlist}
\item[(i)] if LLN-$P${-ex} holds, and then
\begin{eqnarray*}
L(\mu,t,\lambda)&= &
\frac{|\phi'(t)|}{\mu} \biggl(\frac{e^{-\lambda(1-\phi
(t))}-e^{-\lambda}}{\lambda
\phi(t)} \biggr) \\
&= & \frac{|\phi'(t)|}{\mu} \int_0^1 e^{-\lambda\phi(t)x}\, dx
\exp\bigl\{-\lambda\bigl(1-\phi(t)\bigr)\bigr\};
\end{eqnarray*}
\item[(ii)] if LLN-$P^-$ holds, and then
\[
L(\mu,t,\lambda)=\frac{1}{\mu} \int_0^{\infty} x e^{-x t} \exp\biggl\{
-\lambda\int_{x^+}^{\infty} e^{- yt} P(dy)\biggr\} P(dx)\exp\bigl\{-\lambda
\bigl(1-\phi(t)\bigr)\bigr\};
\]
\item[(iii)] if LLN-$P^+$ holds, and then
\[
L(\mu,t,\lambda)=\frac{1}{\mu} \int_0^{\infty} x e^{-x t} \exp\biggl\{
-\lambda\int_0^{x} e^{- yt} P(dy)\biggr\} P(dx)\exp\bigl\{-\lambda\bigl(1-\phi
(t)\bigr)\bigr\}.
\]
\end{longlist}
\end{Th}

The proofs of Theorems \ref{theo:principaleq} and
\ref{theo:principalneq} are deferred to the next section. They
will rely on what we call a law of large numbers for random
partitions of the interval. Let us now deduce the law of the
limiting transient search-cost under the previous sets of
hypotheses.

\begin{coro}\label{coro:limitcaseexchange}
If LLN-$P${-ex} holds, for each $t>0$ we have
\[
\frac{S^{(n)}(n \mu t)}{n} \rightarrow^{d} S(t),
\]
where $S(t)$ satisfies the relation in distribution
%
\begin{equation}\label{decompSt}
S(t)=^{(d)}S^{(\infty)}\1_{\{S^{(\infty)}\leq1-\phi(t)\}}+
U\1_{\{S^{(\infty)}>1-\phi(t)\}},\end{equation}
with $S^{(\infty)}$ defined in Theorem \ref{javiera} and
$U$ a uniform random variable in $[1-\phi(t), 1]$
independent of $S^{(\infty)}$.
Moreover, when $n\to\infty$ we have
\[
\mathbb{P} \bigl(S^{(n)}_{\mathbf{o}}(n \mu t)>0
\bigr)\longrightarrow
\mathbb{P} \bigl(S^{(\infty)}> 1-\phi(t) \bigr)=\frac{|\phi
'(t)|}{\mu}.
\]
Finally, the random variable $S(t)$ has density
\[
f_{S(t)}(x)=f_{S^{(\infty)}}(x)\1_{[0,1-\phi(t)]}+\frac{|\phi
'(t)|}{\mu
\phi(t)}\1_{[1-\phi(t),1]}
\]
and, with $\|\cdot\|_{\mathrm{TV}}$ denoting the total variation distance,
we have
\begin{eqnarray*}
\big\| \operatorname{law}(S(t))-\operatorname{law}\bigl(S^{(\infty)}\bigr)\big\|_{\mathrm{TV}}
&= &
\int_{1-\phi(t)}^{1-\mathbf{p}_0}
 \bigg|f_{S^{(\infty)}}(x)+\frac{|\phi'(t)|}{\mu\phi(t)}  \bigg|\, dx
+\frac{ |\phi'(t)|\mathbf{p}_0}{\mu\phi(t)} \\
&\leq&  2
\frac{|\phi'(t)|}{\mu}.
\end{eqnarray*}
\end{coro}
\begin{pf}
The Laplace transform of $S^{(\infty)}$ is
\[
\mathbb{E}\bigl(\exp\bigl\{-\lambda S^{(\infty)}\bigr\}\bigr)=\frac{1}{\mu}
\int_0^{\infty}\phi''(u)\exp\bigl\{-\lambda\bigl(1-\phi(u)\bigr)\bigr\}\, du \]
(see \cite{BarreraHuilletParoissinMtF2} or Remark
\ref{recovstat}). Now, from Theorem~\ref{theo:principaleq} we
have
\[
\lim_{n\to\infty}A_n(n\mu t,\lambda/n)= \frac{1}{\mu} \int_0^t
\phi''(u) \exp \bigl\{-\lambda\bigl(1-\phi(u)\bigr)  \bigr\}\, du .
\]
Taking $\lambda=0$ and using Remark \ref{SC}, we obtain
\begin{eqnarray*}
\lim_{n\rightarrow\infty} P \bigl(S^{(n)}(n \mu t)\leq|R_{n \mu
t}| \bigr)
&=& \lim_{n\rightarrow\infty}A_n(n\mu
t,0)\\
&=&
1-\frac{(-\phi'(t))}{\mu}\\
&=&
 \mathbb{P}\bigl(S^{(\infty)}\leq1-\phi(t)\bigr).
\end{eqnarray*}

On the other hand, since the Laplace transform of
$\frac{S^{(n)}(n \mu t)}{n}$ conditional on the event $S^{(n)}(n
\mu t)\leq|R_{n \mu t}|$ is given by
\[
\mathbb{E} \biggl(\exp \biggl\{-\lambda\frac{S^{(n)}(n\mu
t)}{n} \biggr\} \Big| S^{(n)}(n \mu t)\leq|R_{n \mu t}|  \biggr) =
\frac{A_n(n\mu t,\lambda/n)}{A_n(n\mu t,0)},\]
we obtain
\begin{eqnarray*}
&&\lim_{n\rightarrow\infty} \mathbb{E} \biggl(\exp \biggl\{-\lambda
\frac{S^{(n)}(n\mu t)}{n} \biggr\} \Big|S^{(n)}(n \mu
t)\leq|R_{n \mu t}|  \biggr)\\
&&\qquad =
\frac{1}{\mu+\phi'(t)} \int_0^t \phi''(u) \exp
 \bigl\{-\lambda\bigl(1-\phi(u)\bigr)  \bigr\}\, du \\
&&\qquad = \mathbb{E} \bigl(\exp\bigl\{-\lambda S^{(\infty)}\bigr\} | S^{(\infty
)} \leq
1-\phi(t) \bigr).
\end{eqnarray*}
Concerning the limiting behavior of $B_n$, we get in a similar way
that
\begin{eqnarray*}
\lim_{n\rightarrow\infty} \mathbb{P} \bigl(S^{(n)}(n \mu t)> |R_{n
\mu
t}| \bigr)&=&-\frac{\phi'(t)}{\mu}\\
&=&\mathbb{P}\bigl(S^{(\infty)}> 1-\phi(t)\bigr),
\end{eqnarray*}
and
\begin{eqnarray*}
\lim_{n\rightarrow\infty} \mathbb{E} \biggl(\exp \biggl\{-\lambda
\frac{S^{(n)}(n\mu t)}{n} \biggr\} \Big|S^{(n)}(n \mu t)> |R_{n
\mu t}|  \biggr)
&= & \lim_{n\rightarrow\infty} \frac{B_n(n\mu
t,\lambda/n)}{B_n(n\mu t,0)}\\
&= & \biggl(\frac{e^{-\lambda(1-\phi(t))}-e^{-\lambda}}{\lambda
\phi(t)} \biggr)\\
&= &
\mathbb{E}(\exp\{-\lambda U\}).
\end{eqnarray*}

Combining the previous limits yields
\begin{eqnarray*}
&&\lim_{n\rightarrow\infty} \mathbb{E} \biggl(\exp \biggl\{-\lambda
\frac{S^{(n)}(n
\mu t)}{n} \biggr\} \biggr)\\
&&\qquad=\lim_{n\rightarrow\infty} A_n(n\mu
t,\lambda/n) + B_n(n\mu t,\lambda/ n)\\
&&\qquad= \mathbb{E} \bigl(\exp\bigl\{-\lambda S^{(\infty)}\bigr\} \1_{\{S^{(\infty
)} \leq
1-\phi(t)\}} \bigr)\\
&&\qquad\quad{} +\mathbb{E}\bigl(\exp\bigl\{-\lambda U\bigr\}\bigr)P\bigl(S^{(\infty)}
> 1-\phi(t)\bigr) \\
&&\qquad= \mathbb{E} \bigl(\exp\bigl\{-\lambda S^{(\infty)}\bigr\} \1_{\{S^{(\infty
)} \leq
1-\phi(t)\}} \bigr)\\
&&\qquad\quad{} +\mathbb{E} \bigl(\exp\{-\lambda U\} \1_{\{S^{(\infty)} >
1-\phi(t)\}}  \bigr)\\
&&\qquad=\mathbb{E} \bigl(\exp-\lambda\bigl\{ S^{(\infty)} \1_{\{S^{(\infty)}
\leq
1-\phi(t)\}} +U \1_{\{S^{(\infty)} > 1-\phi(t)\}} \bigr\} \bigr).
\end{eqnarray*}
From the latter we obtain the density of $S(t)$, and then the
total variation distance to equilibrium. The last asserted
inequality follows from the well-known fact that
$\|\operatorname{law}(X)-\operatorname{law}(Y)\|_{\mathrm{TV}}\leq2 \mathbb{P}\{X\not= Y\}$ for any
coupling of
random variables $(X,Y)$.
\end{pf}

\begin{coro}\label{coro:limitcasemonot}
Define two functions $g_t$ and $\widetilde{g_t}$ by
\[
g_t(y)=\int_0^y e^{-z t}P(dz)\quad \mbox{and}\quad \widetilde
{g_t}(y)=g_t^{-1}(1-y) \qquad\mbox{if LLN-$P^-$ holds},
\]
or by
\[
g_t(y)=\int_{y^+}^{\infty}e^{-z t}P(dz)\quad\mbox{and}\quad \widetilde
{g_t}(y)=(1-g_t)^{-1}(y) \qquad\mbox{if LLN-$P^+$ holds}.
\]
(Here, $g^{-1}$ stands for the generalized inverse of a nondecreasing
right continuous function $g\dvtx\mathbb{R}_+\to\mathbb{R}_+$.) Then,
under LLN-$P^-$ or LLN-$P^+$, for each $t>0$ we have
\[
\frac{S^{(n)}(n \mu t)}{n} \rightarrow^{d} S(t),
\]
where $S(t)$ has the density
%
\begin{equation}\label{densSt'}
f_{S(t)}(x)=\1_{[0,1-\phi(t)]}(x)f_{S^{(\infty)}}(x) +
\1_{[1-\phi(t),1]}(x)\frac{1}{\mu} \widetilde{g_t}(x).
\end{equation}
Moreover, we have
\[
\mathbb{P} \bigl(S^{(n)}_{\mathbf{o}}(n \mu t)>0
\bigr)\longrightarrow
\mathbb{P} \bigl(S(t)> 1-\phi(t) \bigr)= \mathbb{P}
\bigl(S^{(\infty)}>
1-\phi(t) \bigr)=\frac{|\phi'(t)|}{\mu},
\]
when $n\to\infty$, and for all $t\geq0$, $ \|
\operatorname{law}(S(t))-\operatorname{law}(S^{(\infty)})\|_{\mathrm{TV}}\leq2 \frac{|\phi'(t)|}{\mu}.$
\end{coro}

\begin{pf}
If LLN-$P^-$ holds, the result follows by using Theorem
\ref{theo:principalneq} and making the change of variable
$z=1-g_t(x)$ to obtain
\begin{eqnarray*} L(\mu,t,\lambda)&= &\frac{1}{\mu} \int_0^{\infty} x e^{-x
t} \exp\biggl\{
-\lambda\int_{x^+}^{\infty} e^{- yt} P(dy)\biggr\} P(dx)\exp\bigl\{-\lambda
\bigl(1-\phi(t)\bigr)\bigr\} \\
&= & \frac{1}{\mu}\int_{1-\phi(t)}^1 \exp\{-\lambda z\} g_t^{-1}(1-z)\,dz.
\end{eqnarray*}
The remaining case is similar.
\end{pf}

\begin{rem}\label{originaltime}
If $\tilde{S}^{(n)}(t)$ denotes the search cost of the MtF
process in the original time-scale (see Section~\ref{sec2}), it is clear
that the previous results are equivalently stated replacing
$S^{(n)}(n\mu t)$ by
\[
\tilde{S}^{(n)}\bigl(n\mu t/w^{(n)}\bigr),
\]
where $w^{(n)}=\sum_{j=1}^n
w^{(n)}_j$.
\end{rem}

\subsection{Examples and applications}

In what follows, we give the limiting distribution of the
transient search-cost for examples of random or deterministic
request probabilities. The first ones are examples where explicit
computations can be easily done.
\begin{longlist}
\item[(1)] Let $w_i\sim \operatorname{Bernoulli}(p)$, then%
\[
f_{S(t)}(x) = \frac{1}{p}
\1_{[0,p(1-e^{-t}))}(x)+\frac{e^{-t}}{1-p+pe^{-t}}
\1_{[p(1-e^{-t}),1]}(x).
\]
\item[(2)] Let $w_i \sim \operatorname{Gamma}(1,\alpha)$, then%
\[
f_{S(t)}(x) =  \biggl(1+\frac{1}{\alpha} \biggr) (1-x)^{1/\alpha}
\1_{[0,u(t))}(x)+(1+t)^{-1}\1_{[u(t),1]}(x),
\]
with $u(t)=1-(1+t)^{-\alpha}$.
\item[(3)] If $w_i \sim \operatorname{Geometric}(p)$, then%
\[
f_{S(t)}(x)= \frac{2(1-x)-p}{1-p} \1_{[0,u(t))}(x)+ \frac{p
e^{-t}}{1-(1-p)e^{-t}}\1_{[u(t),1]}(x) ,
\]
where $u(t)= \frac{(1-p)(1-e^{-t})}{p+(1-p)(1-e^{-t})}$.

\item[(4)] If $w_i = 1$ or equivalently, $w_i\sim\delta_1$,
we get $f_{S(t)}(x) =1$ (using any of LLN-$P${-ex}, LLN-$P^+$ or
LLN-$P^-$). That is, the limiting search cost is uniform for all
$t\geq0$.
\end{longlist}

The stationary distributions associated with examples~(5)(i) and
(6)(i) below were first studied in Fill \cite{FillMtF96}, whereas
the stationary regimes of examples~(5)(ii) and
(6)(ii) were considered
by Barrera, Huillet and Paroissin in \cite{BarreraHuilletParoissinMtF2}.
The description of the stationary behavior of example~(5)(i) is
also included in Theorem~2 of Jelenkovi\'{c} and Radovanovi\'{c}
\cite{JR}, in the more general context of the
Persistent-Access-Caching (PAC) algorithm introduced therein (see
the more detailed discussion below on the PAC algorithm and the
Last-Recently-Used rule).

\begin{longlist}
\item[(5)]Let $\alpha\in(-1,0)$ and define
\begin{eqnarray*}
P_{\alpha}(dx)
&=& -\frac{1}{\alpha} x^{1/\alpha-1}\1_{[1,\infty
)}(x)\,dx \qquad\mbox{(Pareto law)},\\
\phi(s)
&=&-\frac{1}{\alpha}\int_1^{\infty} e^{-xs} x^{1/\alpha-1}\,
dx,\\
g_t(y)
&=&-\frac{1}{\alpha}\int_1^{y} e^{-xt} x^{1/\alpha-1}\,dx.
\end{eqnarray*}
\begin{longlist}\vspace*{-4pt}
\item[(i)] If $w_i = i^{\alpha}$, we have using~\eqref{Pexamplec} that
\begin{eqnarray*}
f_{S(t)}(x)
&=&
-(\alpha+1)\frac{\phi''(\phi^{-1}(1-x))}{\phi'(\phi^{-1}(1-x))}\1
_{[0,1-\phi(t))}(x)\\
&&{}+
(\alpha+1)g_t^{-1}(1-x)\1_{[1-\phi(t),1]}(x).
\end{eqnarray*}
\item[(ii)] If $w_i$ are i.i.d. with law $P_{\alpha}$, then
\begin{eqnarray*}
f_{S(t)}(x)&=&-(\alpha+1)\frac{\phi''(\phi^{-1}(1-x))}{\phi'(\phi
^{-1}(1-x))}\1_{[0,1-\phi(t))}(x)\\
&&{}+
(\alpha+1)\frac{|\phi'(t)|}{\phi(t)}\1_{[1-\phi(t),1]}(x).
\end{eqnarray*}
\end{longlist}
\item[(6)]Let $\alpha> 0$ and set now
\begin{eqnarray*}
\phi(s)
&=&\frac{1}{\alpha}\int_0^1 e^{-xs} x^{1/\alpha-1}\,
dx,\\
g_t(y)
&=&\frac{1}{\alpha}\int_y^1 e^{-xt} x^{1/\alpha-1}\,dx.
\end{eqnarray*}
\begin{longlist}
\item[(i)] If $w_i = i^{\alpha}$, we have by~\eqref{Pexamplec}
that
\begin{eqnarray*}
f_{S(t)}(x)
&=&-(\alpha+1)\frac{\phi''(\phi^{-1}(1-x))}{\phi'(\phi^{-1}(1-x))}\1
_{[0,1-\phi(t))}(x)\\
&&{}+
(\alpha+1)(1-g_t)^{-1}(x)\1_{[1-\phi(t),1]}(x).
\end{eqnarray*}
\item[(ii)] If $w_i$ are i.i.d. with law $\operatorname{Beta}(1,1/\alpha)$, then
\begin{eqnarray*}
f_{S(t)}(x)
&=&-(\alpha+1)\frac{\phi''(\phi^{-1}(1-x))}{\phi'(\phi
^{-1}(1-x))}\1_{[0,1-\phi(t))}(x)\\
&&{}+
(1+\alpha)\frac{|\phi'(t)|}{\phi(t)}\1_{[1-\phi(t),1]}(x).
\end{eqnarray*}
\end{longlist}
\end{longlist}

It was remarked in \cite{BarreraHuilletParoissinMtF2} that
example~(5)(i) shares the same stationary distribution as example~(5)(ii), and example~(6)(i)
the same as that of (6)(ii). We observe
here that the equilibrium parts of the transient search-costs of
examples~(5)(i) and~(5)(ii) coincide as well, as happens also with
examples~(6)(i) and (6)(ii). In turn, their out of equilibrium
transient search-cost are different. In Section~\ref{sec5} we discuss
and explain these facts in the light of the new techniques that
will be shortly introduced.

To motivate our last example, we recall that Jelenkovi\'{c}
\cite{Jelen} considered a continuum (thus infinite) list of
objects representing an ``efficient static'' or popularity
decreasing arrangement of objects. More precisely, a probability
measure $Q$ on~$\mathbb{R}_+$ with decreasing density $q$ is used therein
to specify the probability $q(x)\,dx$ that an object lying at
position $x\in\mathbb{R}_+$ is requested. The stationary search-cost was
studied by approximating the continuum list by discrete albeit
countably infinite lists $(Q^n)_{n\in\mathbb{N}}$, (a ``fluid limit'').
That is, for each $n$, object at position $i/n$, $i\in\mathbb{N}$, is
requested with probability $Q^n(i)=Q((i+1)/n)-Q(i/n)$. In this
case the transient search cost $S(t)$ can be obtained by our
approach when the continuum list has finite length.

\begin{longlist}
\item[(7)] Let $Q(dx)=q(x)\,dx$ be supported in $[0,c]$ and
$w^{(n)}_i,i=1,\ldots,n$, be defined as in {(d)}
above. If $q$ is a continuous decreasing probability density,
then LLN-$P^-$ holds and we get
\[
f_{S(t)}(x) =-c
\frac{\phi''(\phi^{-1}(1-x))}{\phi'(\phi^{-1}(1-x))}\1_{[0,1-\phi
(t))}(x)\\+
cg_t^{-1}(1-x)\1_{[1-\phi(t),1]}(x),
\]
where
\begin{eqnarray*}
\phi(s)
&=&\frac{1}{c}\int_0^c e^{-q(x)s} \,dx \quad\mbox{and}\\
g_t(y)
&=&\frac{1}{c}\int^{y}_0 e^{-q(x)t} \,dx.
\end{eqnarray*}
Observe that in this case
$\mu=1/c$ since $P(dx)=(\frac{1}{c}\1_{x\in[0,c]}\,dx)\circ q^{-1}$, so that the Laplace transform
computed in Theorem~\ref{theo:principaleq} reads
\[
\int_{0}^{t} \biggl(\int_{0}^c
q^2(u)e^{-q(u)s}\,du \biggr)\exp \biggl(-\lambda\biggl(1-\frac{1}{c}\int
_{0}^ce^{-q(u)s}\,du\biggr) \biggr)\,ds.
\]
The limit of this expression when $t\to\infty$ is exactly
formula 4.1 in~\cite{Jelen} evaluated in $s=\lambda/c$.
\end{longlist}

We notice that the fluid limit approximation of Jelenkovi\'{c}
\cite{Jelen} also corresponds to a law of large numbers
asymptotic, in the sense that the (infinite) empirical measures
$\frac{1}{n}\sum_{i\in\mathbb{N}} \delta_{i/n}$ approach $dx$ in
$\mathbb{R}_+$
as the space scale $1/n$ goes to $0$. However, in our case the
search-cost is defined in terms of the relative position in a
finite list, whereas in~\cite{Jelen} it is understood as the
absolute position in a possibly infinite list. (The reader
familiar with particle systems will recognize a similar difference
between the hydrodynamic and mean field limit formalisms; it is from
the latter that we have borrowed the law of large numbers
terminology; see next section.) Although we can ``simulate'' the
fluid limit for compactly supported measures $Q$ [as example~(7)
shows] the general case is not tractable with our techniques (see
the discussion in Section~\ref{sec5}).

To finish the discussion on related works, we remark that our
main results also describe the transient behavior of particular
instances of the Least-Recently-Used (LRU) caching rule, which
dynamically selects a collection of frequently accessed documents
and stores them in a low cost access place. Indeed, the
probability that at time $n\mu t$ the requested document is not
found among the $\delta n$ selected ones (and a fault occurs)
corresponds to the probability $\mathbb{P}(S^{(n)}(n\mu t)>\delta n)$
that the search-cost in the MtF scheme is bigger than $\delta n$
(for a more detailed discussion of this relation we refer to
\cite{Jelen}). Consequently, from Corollary~\ref{coro:limitcaseexchange} we can,
for instance, compute the
transient asymptotic fault probability under assumption
LLN-$P$-{ex}:
\begin{eqnarray}\label{eq:LRU}
\mathbb{P}\bigl(S(t)>\delta\bigr) =
\cases{
\dfrac{|\phi'(\eta_{\delta})|}{\mu}, &\quad if $\eta_{\delta} < t$,\vspace*{2pt}\cr
\dfrac{1-\delta}{\mu}\frac{\phi'(t)}{\phi(t)},&\quad if $\eta_{\delta} \geq t$,
}
\end{eqnarray}
with $\eta_{\delta}= \phi^{-1}(1-\delta)$. Under LLN-$P^+$ or
LLN-$P^-$, thanks to Corollary~\ref{coro:limitcasemonot}, the same value is
obtained for the case
$\eta_{\delta}<t$, and an integral expression in terms of~$\tilde
{g}_t$ (which can be written
explicitly)
in the case $\eta_{\delta}\geq t$.

We remark that the PAC algorithm introduced in \cite{JR}
generalizes the LRU rule by updating the list in a similar way,
but only if the requested item at time $t$ has already been
requested $k-1$ times in the time interval $((t-\beta)\vee0,t)$.
Thus, by taking $t=\infty$ in formula~\eqref{eq:LRU} (see Remark~\ref{recovstat})
we obtain a generalization of the stationary
result of Theorem 2 of \cite{JR} in the case $k=1,\beta>0$ [the
latter corresponding to the particular $\phi$ given by the Pareto
law of example~(5)(i)]. Moreover, for the case $k=1$ and
$w_i=i^{\alpha}$, $\alpha\in(-1,0)$ of Theorem 2 of \cite{JR},
we obtain the transient asymptotic fault
probability. This is given by
\begin{eqnarray*}
&&\mathbb{P}\bigl(S(t)>\delta\bigr)\\
&&\qquad=\cases{
-\dfrac{\alpha+1}{\alpha}\eta_{\delta}^{-(1+1/\alpha)}\Gamma
 (1+1/\alpha,\eta_{\delta} ),  & \quad if $\eta
_{\delta} < t$,\vspace*{2pt}\cr
-\dfrac{\alpha+1}{\alpha}t^{-(1+1/\alpha)} [\Gamma
(1+1/\alpha,t  )-\Gamma (1+1/\alpha,t\varepsilon_{\delta
,t} ) ],&\quad  if $\eta_{\delta} \geq t$,
}
\end{eqnarray*}
where $\Gamma(z,y):=\int_y^{\infty} x^{z-1}e^{-x} \,dx$ is the
incomplete Gamma function, and
$\varepsilon_{\delta,t}:=g_t^{-1}(1-\delta)$ with $g_t$ as in example~(5)(i).

\subsection{Stochastic order relations}

In the remainder of this section we shall establish some
stochastic order relations between the three situations LLN-$P${-ex},
LLN-$P^+$ and LLN-$P^-$. Recall that
given two real valued random variable $X$ and $Y$, we say that
\textit{$X$ is stochastically smaller than $Y$}, if for all $z\in
\mathbb{R}$, one has $\mathbb{P}(X\leq z)\geq\mathbb{P}(Y\leq z)$.
This is written $X
\preceq Y$.

Notice now that the three assumptions can be seen as a priori
information of different type about the initial positions of
objects in the list. More precisely, LLN-$P${-ex} can be
read as having no a priori knowledge at all, whereas
LLN-$P^-$ can be interpreted as the relative order of
popularities being known, and objects being placed at time $0$
in decreasing order (intuitively, this is an efficient statical
ordering). Accordingly, assumption LLN-$P^+$ can be
interpreted as the least efficient order at time $0$, if the
relative order of popularities is known. In this direction, Fill
and Holst proved in Corollary 4.2 of \cite{FillHolst96} that for a
given finite request probability vector, the transient search-cost
is stochastically larger than that of the same vector rearranged
in decreasing order, and smaller than when it is arranged in
increasing order. We shall prove that similar stochastic order
relations hold in the large numbers limit, by using the explicit
expressions for $f_{S(t)}$ we have already found.

\begin{coro}
Let $S^{ex}(t)$, $S^+(t)$ and $S^-(t)$ denote the limiting
transient search-cost $S(t)$ respectively under the assumptions,
LLN-$P${-ex}, LLN-$P^+$ and  LLN-$P^-$. Then,
we have
\[
S^-(t)\preceq S^{ex}(t)\preceq S^+(t).
\]
\end{coro}
\begin{pf}
From Corollaries~\ref{coro:limitcaseexchange} and~\ref{coro:limitcasemonot}, we just need to prove that
\begin{eqnarray*}
\mathbb{P}\{1-\phi(t)\leq S^-(t)\leq x\}&\geq&\mathbb{P}\{1-\phi
(t)\leq
S^{ex}(t)\leq x\}\\
&\geq&\mathbb{P}\{ 1-\phi(t)\leq S^+(t) \leq
x\}
\end{eqnarray*}
for all $x\in[1-\phi(t),1]$. The first
inequality is equivalent to
\[
\int_{g_t^{-1}(1-x)}^{\infty}z e^{-zt}P(dz)\geq
\frac{|\phi'(t)|}{\phi(t)}\bigl(x-1+\phi(t)\bigr)
\]
for all $x\in[1-\phi(t),1]$, where $g_t(y)=\int_0^y e^{-z
t}P(dz)$. This will follow if we can prove that
\[
\frac{\int_{y+}^{\infty} z e^{-zt}P(dz)}{|\phi'(t)|}\geq
\frac{\int_{y+}^{\infty} e^{-zt}P(dz)}{\phi(t)}
\]
for all $y\geq0$ or, equivalently, that
%
\begin{equation}\label{ineqstoc}
\frac{\int_0^{y} z e^{-zt}P(dz)}{|\phi'(t)|}\leq\frac{\int_0^{y}
e^{-zt}P(dz)}{\phi(t)}.
\end{equation}
Observe that both sides have the same points of discontinuity, as
functions of $y$. Therefore, by suitably approximating $P$, we may
assume that $P(dz)$ has a continuous density $f(z)$ which is
strictly positive. Write $a(y)=\int_0^{y} z e^{-zt}f(z)\,dz$ and
$b(y)=\int_0^{y} e^{-zt}f(z)\,dz$. We need to check that
\[
h(y):= \frac{a(y)}{a(\infty)}- \frac{b(y)}{b(\infty)}\leq0.
\]
Since $h$ is differentiable and $h(0)=h(\infty)=0$, it is enough
to prove that $h$ has a unique critical point $y_0$ and that
$h(y_0)\leq0$. By the assumption on $P$, the condition
$h'(y_0)=0$ is satisfied if and only if
$y_0=\frac{a(\infty)}{b(\infty)}$. But then, $h(y_0)\leq0$ is the
same as
\[
\frac{\int_0^{a(\infty)/b(\infty)} z
e^{-zt}f(z)\,dz}{a(\infty)}\leq
\frac{\int_0^{a(\infty)/b(\infty)}
e^{-zt}f(z)\,dz}{b(\infty)},
\]
which is trivially true. We conclude that $\mathbb{P}\{1-\phi(t)\leq
S^-(t)\leq x\}\geq\mathbb{P}\{1-\phi(t)\leq S^{ex}(t)\leq x\}$ for all
$x\geq0$. The remaining inequality is easily seen to follow also
from~(\ref{ineqstoc}).
\end{pf}

\section{Law of large numbers for random partitions of the interval
and proofs of
Theorems~\protect\ref{theo:principaleq} and \protect\ref{theo:principalneq}}\label{sec4}

The main ingredient in the proofs of Theorems
\ref{theo:principaleq} and~\ref{theo:principalneq} will be what we
call a ``law of large numbers for random partitions of the
interval.'' To illustrate this idea, consider first $(w_i)_{i\in
\mathbb{N}}$ i.i.d. random variables in $\mathbb{R}_+$ of law $P$
with finite
mean $\mu>0$, and the
probability vector $\mathbf{p}^{(n)}=(p_i^{(n)})$ defined by
\[
p_i^{(n)}:=\frac{w_i}{\sum_{j=1}^n w_j},\qquad i=1,\ldots, n.
\]
Then, by the strong law of large numbers, we have
\[
\bigl(n \mu p^{(n)}_1,\ldots,n \mu p^{(n)}_k\bigr)\longrightarrow(w_1,\ldots,w_k)
\]
almost surely when $n\to\infty$.
In particular, any $k\leq n$ fixed coordinates of the vector
$n\mu\mathbf{p}^{(n)}$ become independent as $n$ tends to
infinity, and the limiting law of each of them converges to $P$.
The following result due to H. Tanaka implies that the empirical
measures
\[
\nu^{(n)}:=\frac{1}{n}\sum_{i=1}^n \delta_{n \mu p^{(n)}_i},
\]
converge to $P$, as $n$ goes to infinity.

\begin{prop}\label{propchaos}
For each $n\in\mathbb{N}$, let $X^{(n)}=(X^{(n)}_1,\ldots,X^{(n)}_n)$ be
an exchangeable random vector in $\mathbb{R}^n$ with law $P_n$. Then,
the following assertions are equivalent:
\begin{longlist}
\item[(i)] There exists a probability measure $P$ in $\mathbb{R}$
such that for all $k\in\mathbb{N}$, when $n\to\infty$,
\[
\operatorname{law} \bigl(X^{(n)}_1,\ldots,X^{(n)}_k\bigr)\Longrightarrow P^{\otimes k}.
\]
\item[(ii)]
The random variables $\frac{1}{n}\sum_{i=1}^n \delta_{X_i^{(n)}}$
[taking values in the polish space
$\mathcaligr{P}(\mathbb{R})$] converge in law as $n$ goes to infinity
to a
deterministic limit equal to $P$.
\end{longlist}
\end{prop}

A sequence of probability measures $P_n$ satisfying condition~{(i)} of Proposition~\ref{propchaos} is said to be $P$-chaotic, or
to have the propagation of chaos property with limiting law $P$.
This is a central property in the probabilistic study of mean
field models. For further background on these topics and a proof
of Proposition~\ref{propchaos}, we refer the reader to Sznitman's
course \cite{Szn}.

We now prove that the same conclusion about $\nu^{(n)}$ can be
obtained under a weaker assumption on the vectors
$(w_i^{(n)})_{i=1}^n$, $n\in\mathbb{N}$. Namely, we have:

\begin{Th}[(L.L.N. for random partitions of the interval)]
\label{teoLLN}
Assume that
$(\mathbf{w}^{(n)})_{n\in\mathbb{N}}$ satisfy condition LLN-$P$
and let $(\mathbf{p}^{(n)})_{n\in\mathbb{N}}$ be defined as in~\eqref{randompartition}. Then, the empirical measure
\[
\nu^{(n)}:=\frac{1}{n}\sum_{i=1}^n \delta_{n \mu p^{(n)}_i}
\]
converges in law to the deterministic limit $P$.
\end{Th}

This amounts to say that if LLN-$P$ holds for
$(w_i^{(n)})_{i=1}^n$ and some sequence $(Z_n)$, then it also
holds for $(p_i^{(n)})_{i=1}^n$ and the sequence $(Z_n'):=(n\mu)$
(the convergence of the empirical means of $\nu^{(n)}$ being
trivial).

\begin{pf*}{Proof of Theorem \ref{teoLLN}}
The proof is simple by using the Wasserstein distance
$W_1$ in the space ${\mathcaligr
P}_1(\mathbb{R})$ of Borel probability measures on $\mathbb{R}$ with finite
first moment. Recall that
\[
W_1(m,m')=\inf_Q \int_{\mathbb{R}^2} |x-y| Q(dx,dy),
\]
where the $\inf$ is taken over all Borel probability measures $Q$
on $\mathbb{R}^2$ with first and second marginal laws in ${\mathcaligr
P}_1(\mathbb{R})$ respectively equal to $m$ and $m'$ (i.e., couplings
of~$m$ and $m'$). Then, $W_1$ is a distance inducing the weak
topology, strengthened with the convergence of first-order moments
(see, e.g., \cite{Vil}). Let us define
\[
Q^{(n)}:=\frac{1}{n}\sum_{i=1}^n \delta_{ (n \mu
p^{(n)}_i,Z_nw_i^{(n)} )}
\]
which is a coupling of $\hat{\nu}^{(n)}$ and $\nu^{(n)}$. Then, on
the event $\{\sum_{j=1}^n Z_nw_j^{(n)}>0\}$ we have
\begin{eqnarray*}
\int_{\mathbb{R}^2} |x-y| Q^{(n)}(dx,dy)
&=& \frac{1}{n}\sum_{i=1}^n  \big| n \mu p^{(n)}_i-Z_nw_i^{(n)}\big| \\
&= & \frac{1}{n}\sum_{i=1}^n Z_nw_i^{(n)} \bigg| \frac{n\mu
}{\sum_{j=1}^n Z_nw_j^{(n)}}-1\bigg |\\
&= &  \bigg|\mu- \frac{\sum_{i=1}^n Z_nw_i^{(n)}}{n }\bigg |,
\end{eqnarray*}
from where
\[
W_1\bigl(\nu^{(n)},\hat{\nu}^{(n)}\bigr)\leq \bigg| \frac{\sum_{i=1}^n
Z_nw_i^{(n)}}{n }-\mu\bigg |.
\]
We deduce that
\[
W_1\bigl(\nu^{(n)},P\bigr)\leq \bigg| \frac{\sum_{i=1}^n Z_nw_i^{(n)}}{n
}-\mu \bigg| +W_1\bigl(\hat{\nu}^{(n)},P\bigr).
\]
Now, from LLN-$P$ we have $\frac{\sum_{i=1}^n
Z_nw_i^{(n)}}{n }\to\mu$ and $\hat{\nu}^{(n)}\to P$, both in\break
probability (the second with respect to $W_1$). On the other hand,
we have $\mathbb{P}\{ \sum_{j=1}^n Z_nw_j^{(n)} =0\}\to0$ as $n\to
\infty$, since $\mu>0$. We deduce that $\nu^{(n)}$ converges to
$P$ in probability with respect to $W_1$, and therefore also with
respect to the usual weak topology.
\end{pf*}

\begin{rem}\label{conseqLLN}
Assumption LLN-$P$ together with Theorem \ref{teoLLN}, imply
that for any bounded continuous function $\Psi\dvtx{\mathcaligr P}(\mathbb{R}
_+)\to
\mathbb{R}$, one has
\[
\mathbb{E}(\Psi(\nu^{n}))\to\Psi(P)\qquad\mbox{when }n\to\infty.
\]
We shall systematically rely on this fact in the proof of Theorems~\ref{theo:principaleq} and \ref{theo:principalneq}, in order to
compute the limits of quantities of the form
$\mathbb{E} (f_n(\mathbf{p}^{(n)}) )$, for adequate functions
$f_n\dvtx\mathbb{R}_+^{n}\to\mathbb{R}$, $n\in\mathbb{N}$. Namely, we
will use the fact
that
\[
\mathbb{E} \bigl(f_n\bigl(\mathbf{p}^{(n)}\bigr) \bigr)\to\Psi(P)
\]
whenever
$f_n(\mathbf{p}^{(n)})$ is equal (or asymptotically close enough)
to $\Psi(\nu^{n})$, for some bounded continuous $\Psi$ not
depending on $n$.
\end{rem}

We shall also need the following lemma on the size-biased
picking of probability measures on $\mathbb{R}^+$. Recall that given
$\mathbf{m}\in{\mathcaligr P}(\mathbb{R}_+)$ with $0<\langle
\mathbf{m},x\rangle<\infty$, the \textit{size-biased picking of}
$\mathbf{m}$ is the law
\[
\overline{\mathbf{m}}(dy):=\frac{y}{\langle\mathbf{m},x\rangle
}{\mathbf{m}}(dy);
\]
it is obtained from $\mathbf{m}$ as in the \textit{waiting time
paradox} (see, e.g., Feller \cite{Fe}, Chapter~VI).

\begin{lem}\label{contsbp}
Let $(\mathbf{m}_n)$ be a sequence of probability measures on
$\mathbb{R}_+$ with finite means and weakly converging to a probability
measure $\mathbf{m}\not= \delta_0$. Assume moreover that $\langle
\mathbf{m},x\rangle<\infty$ and that $\langle
\mathbf{m}_n,x\rangle\to\langle\mathbf{m},x\rangle$ when $n$
goes to $\infty$. Then, we have
\[
\overline{\mathbf{m}}_n\Longrightarrow\overline{\mathbf{m}}.
\]
\end{lem}
\begin{pf} Since $\langle\overline{\mathbf{m}}_n,1\rangle
\longrightarrow\langle
\overline{\mathbf{m}},1\rangle$ as $n$ goes to $\infty$, it is
enough to prove that $\langle
\overline{\mathbf{m}}_n,f\rangle\longrightarrow\langle
\overline{\mathbf{m}},f\rangle$ for each continuous function $f$
with compact support. Since for such $f$ the function $x f(x)$ is
continuous and bounded, this follows from the assumptions.
\end{pf}

\begin{pf*}{Proof of Theorem \ref{theo:principaleq}}
In what follows, we drop for notational simplicity the superscript
$(n)$ of the popularity $p_i^{(n)}=p_i$.

From Proposition \ref{prop:bodell}, it holds that
\[
A_n(n \mu t,\lambda/n)=\mathbb{E} \Biggl( \int_0^{n\mu t} e^{-u}
\sum_{i=1}^n p_i^2 \Biggl(\prod_{j=1, j\neq i }^n  \bigl(1+(e^{p_j
u}-1)e^{-\lambda/n} \bigr) \Biggr)\, du  \Biggr).
\]
Let $g_{n},h_{n}\dvtx\mathbb{R}^2_+\to\mathbb{R}$ be functions defined
by
\begin{eqnarray*}
g_n(u,x)&=& \frac{x^2e^{-xu}}{1-(1-e^{-xu})(1-e^{-\lambda/n})},\\
h_n(u,x)&=& n\log \bigl(1-(1-e^{-xu})(1-e^{-\lambda/n}) \bigr).
\end{eqnarray*}
Making the right change of variable we can write $A_n(n \mu
t,\lambda/n)$ as
\begin{eqnarray*}
A_n(n \mu t,\lambda/n)
& = &
\frac{1}{\mu}\mathbb{E} \Biggl( \int_0^t
\frac{1}{n}\sum_{i=1}^n g_n(u,n \mu p_i) \exp \Biggl(
\frac{1}{n}\sum_{i=1}^n h_n(u,n \mu p_i) \Biggr)\, du  \Biggr) \\
& =&
 \frac{1}{\mu}\mathbb{E} \Biggl( \int_0^t \big\langle\nu^{(n)},
g_n(u,\cdot) \big\rangle\exp \bigl( \big\langle\nu^{(n)} ,
h_n(u,\cdot)\big\rangle \bigr)\, du  \Biggr).
\end{eqnarray*}
Now define $g(u,x)= x^2e^{-x u}$ and $ h(u,x)= -(1-e^{-x u})\lambda$.
Then, if
\[
\widetilde{A}_n(n \mu t,\lambda/n):= \frac{1}{\mu}\mathbb{E}
 \biggl( \int_0^t \big\langle\nu^{(n)}, g(u,\cdot) \big\rangle\exp \bigl(
\big\langle\nu^{(n)} , h(u,\cdot)\big\rangle \bigr)\, du  \biggr),
\]
we see
that
\[
|\widetilde{A}_n(n \mu t,\lambda/n)-A_n(n \mu t,\lambda/n)| \leq
I_1(t,\lambda)+I_2(t,\lambda),
\]
with $I_1(t,\lambda)$ and $I_2(t,\lambda)$ defined by
\begin{eqnarray*}
I_1(t,\lambda)
&=& \frac{1}{\mu}\mathbb{E} \biggl( \int_0^t \big\langle
\nu^{(n)},
|g_n(u,\cdot)-g(u,\cdot)|\big\rangle\exp \bigl(
\big\langle\nu^{(n)} , h_n(u,\cdot)\big\rangle \bigr)\, du  \biggr),\\
I_2(t,\lambda)
&=& \frac{1}{\mu}\mathbb{E} \biggl(\int_0^t\big\langle
\nu^{(n)},
g(u,\cdot) \big\rangle \big| \exp \bigl( \big\langle\nu^{(n)} ,
h_n(u,\cdot)-h(u,\cdot)\big\rangle \bigr)-1 \big|\, du  \biggr).
\end{eqnarray*}
On the other hand,
we have the following estimates for $n$ large enough:
%
\begin{eqnarray}\label{estimcabron1}
|g_n(u,x)-g(u,x)|
&\leq& g(u,x)
\frac{1-e^{-\lambda/n}}{1-(1-e^{-ux})(1-e^{-\lambda/n})}\nonumber\\[-8pt]\\[-8pt]
 &\leq&
2g(u,x)\lambda/n\nonumber
\end{eqnarray}
(we use the bound $\frac{1-e^{-\alpha}}{1-c(1-e^{-\alpha})}\leq
e^{\alpha}-1$ for $c\in[0,1], \alpha\geq0$), and
%
\begin{eqnarray}\label{estimcabron2}\hspace*{20pt}
|h_n(u,x)-h(u,x)|
&\leq& 2\lambda
 \biggl\{ \biggl(\frac{\log(1-(1-e^{-xu})(1-e^{-\lambda
/n}))}{(1-e^{-xu})(1-e^{-\lambda/n})}
+1 \biggr)+ \lambda/n \biggr\}\nonumber\\[-8pt]\\[-8pt]\hspace*{20pt}
&\leq& \frac{8\lambda^2}{n}.\nonumber
\end{eqnarray}
In the last line, we have used the bound
\[
 \bigg|\frac{\log(1-c(1-e^{-\alpha}))}{c(1-e^{-\alpha})}+1
\bigg|\leq
2(1-e^{-\alpha})
\]
for all $c\in[0,1]$ and $(1-e^{-\alpha})\leq
1/2$.

Estimates (\ref{estimcabron1}) and (\ref{estimcabron2}) imply that
for large enough $n$, we have
\[
I_1(t,\lambda) \leq\frac{2}{\mu}\mathbb{E} \biggl( \int_0^t
\big\langle\nu^{(n)},
g(u,\cdot)\big\rangle\, du \biggr)\lambda/n\]
and
\begin{eqnarray*}
I_2(t,\lambda)
&\leq&
 \frac{1}{\mu}\mathbb{E} \biggl( \int_0^t
\big\langle\nu^{(n)}, g(u,\cdot)\big\rangle\big\langle\nu^{(n)} ,
|h_n(u,\cdot)-h(u,\cdot)|\big\rangle\, du
 \biggr) \\
&\leq&
 \frac{8}{\mu}\mathbb{E} \biggl( \int_0^t \big\langle\nu^{(n)},
g(u,\cdot)\big\rangle\, du  \biggr)\frac{\lambda^2}{n}.
\end{eqnarray*}
Since $\int_0^t \langle\nu^{(n)}, g(u,\cdot)\rangle\,
du=\int_{\mathbb{R}_+} x(1-e^{-xt})\nu^{(n)}(dx)\leq\mu$ by Fubini's
theorem, we get from the previous estimates that
\[
|\widetilde{A}_n(n \mu t,\lambda/n)-A_n(n \mu t,\lambda/n)| \leq
\frac{C}{n}
\]
for all $n$ large enough. Consequently, we just need to prove that
%
\begin{equation}\label{limtildeA}
\lim_{n\to\infty} \widetilde{A}_n(t,s)= \frac{1}{\mu}\mathbb
{E} \biggl(
\int_0^t \langle P, g(u,\cdot) \rangle\exp\langle P ,
h(u,\cdot)\rangle\, du  \biggr).
\end{equation}

Let us set
\begin{eqnarray*}
\triangle(t,n)&:=&   \bigg| \mathbb{E} \biggl( \int_0^t \big\langle
\nu^{(n)}, g(u,\cdot) \big\rangle\exp \bigl( \big\langle\nu^{(n)} ,
h(u,\cdot)\big\rangle
 \bigr)\, du  \biggr) \\
 &&\hspace*{31.5pt}{}-\int_0^t \langle P, g(u,\cdot) \rangle\exp
\langle P, h(u,\cdot)\rangle\, du  \bigg|.\end{eqnarray*}

For each $\delta>0$, since $h(u,x)\leq0$ we
have the estimate
\begin{eqnarray*}
\triangle(t,n)
&\leq&
  \bigg| \mathbb{E} \biggl( \int_{\delta}^t
\big\langle\nu^{(n)}, g(u,\cdot)\big \rangle\exp\big\langle\nu^{(n)} ,
h(u,\cdot)\big\rangle\, du  \biggr) \\
&&\hspace*{25pt}{}  -\int_{\delta}^t \langle P, g(u,\cdot) \rangle\exp
\langle P, h(u,\cdot)\rangle\, du  \bigg|\\
&&{}  + \int_0^{\delta} \mathbb{E}\big\langle\nu^{(n)}, g(u,\cdot)\big
\rangle\, du + \int_0^{\delta}
\langle P, g(u,\cdot) \rangle\, du .
\end{eqnarray*}

Observe that for each $u>0$ the functions $g(u,\cdot)$ and
$h(u,\cdot)$ are continuous and bounded. Moreover, for each
$\delta>0$, the restriction of $g$ to $[\delta,\infty] $ is
uniformly bounded. Thus, by using dominated convergence, the
mapping
\[
\nu\mapsto F(\nu):=\int_{\delta}^t \langle\nu, g(u,\cdot)
\rangle\exp\langle\nu, h(u,\cdot)\rangle\, du
\]
is seen to be continuous and bounded on $\in{\mathcaligr P}(\mathbb
{R}_+) $.
Thanks to LLN-$P$ and Theorem~\ref{teoLLN}, we deduce that
\[
\mathbb{E}\bigl(F\bigl(\nu^{(n)}\bigr)\bigr)\to F(P)\qquad \mbox{when $n$ goes to
$\infty$}
\]
and, consequently, we get that for any $\delta>0$
%
\begin{equation}\label{limsuptriangle}
\limsup_{n\to\infty} \triangle(t,n) \leq\sup_{n\in\mathbb{N}}
\int_0^{\delta} \mathbb{E}\big\langle\nu^{(n)}, g(u,\cdot)\big \rangle \,du +
\int_0^{\delta} \langle P, g(u,\cdot) \rangle\, du.
\end{equation}

In order to prove (\ref{limtildeA}) it is therefore enough to
establish that the two terms on the r.h.s. of inequality~(\ref{limsuptriangle}) go to $0$ with $\delta$. Notice that the
second term is equal to
\begin{eqnarray*}
\int_{\mathbb{R}_+}
 \biggl( \int_0^{\delta} x^2 e^{-xu}\, du  \biggr) P(dx)
&=& \int_{\mathbb{R}_+} x P(dx) - \int_{\mathbb{R}_+} x e^{-x
\delta} P(dx) \\
&=&  \mu\bigl(\overline{\phi} (0)- \overline{\phi}(\delta)\bigr),
\end{eqnarray*}
where $\overline{\phi}(s):=\frac{1}{\mu}\int_{\mathbb{R}_+}x e^{-s
x}P(dx)$ is the Laplace transform of the size-biased picking of
$P$. Thus, that term goes to $0$ with $\delta$ by continuity of
$\bar{\phi}$.

To tackle the first term on the r.h.s in (\ref{limsuptriangle}),
we consider the intensity measures associated with the random
measures $\nu^{(n)}$. That is, the (deterministic) probability
measures defined for each $n\in\mathbb{N}$ by
\[
\langle\mathbf{m}_n, f \rangle:= \mathbb{E}\big\langle\nu^{(n)},
f\big\rangle .
\]
Notice that $\mathbf{m}_n$ has mean $\mu$ for all $n\in\mathbb{N}$. On
the other hand, if we denote by $\overline{\mathbf{m}}_n$ the
size-biased picking of $\mathbf{m}_n$, we get through
similar computations as before that
\[
\int_0^{\delta} \mathbb{E}\big\langle\nu^{(n)}, g(u,\cdot)\big \rangle\, du
= \mu
\bigl(\overline{\phi}_n (0)- \overline{\phi}_n(\delta)\bigr),
\]
with $\overline{\phi}_n(s):=\frac{1}{\mu}\int_{\mathbb{R}_+}x e^{-s
x}\mathbf{m}_n(dx)$ the Laplace transform of\vspace*{1.5pt}
$\overline{\mathbf{m}}_n$.

Consequently, what we need to prove is that
%
\begin{equation}\label{convlaplaceunif}
\lim_{\delta\to0} \sup_{n \in\mathbb{N}}|\overline{\phi}_n
(\delta)-
\overline{\phi}_n(0)|=0.
\end{equation}
But from LLN-$P$ and Theorem \ref{teoLLN}, for all $f\in
C_b(\mathbb{R})$ we have that
\[
\langle\mathbf{m}_n, f \rangle=\mathbb{E}\big\langle\nu^{(n)},
f\big\rangle
\to\langle P,f \rangle
\]
since the mapping $\nu\mapsto\langle\nu, f\rangle$ is continuous and
bounded. In other words, the sequence $\mathbf{m}_n$ converges
weakly to $P$. With Lemma~\ref{contsbp} we deduce that the
sequence~$\overline{\mathbf{m}}_n$ is weakly convergent, and
therefore, by standard properties of the Laplace transform, the
family of functions $(\overline{\phi}_n)_{n \in\mathbb{N}}$ is
equicontinuous. Clearly, this implies that~(\ref{convlaplaceunif}) holds, and the proof is finished.
\end{pf*}

In the remaining proof we shall use the following result.

\begin{lem}\label{lemma:continvdistrfunc}
Let $F_m$ denote the distribution function of $m\in{\mathcaligr
P}(\mathbb{R})$, and
\[
F_m^{-1}(x):=\inf\{t\geq0\dvtx F_m(t)\geq x\}
\]
be
its generalized inverse. Assume that $m_k\in{\mathcaligr P}(\mathbb{R})$
converges weakly to $m$. Then, $F_{m_k}^{-1}(x)$ converges to
$F_{m}^{-1}(x)$ for $dx$---almost every $x\in[0,1]$.
\end{lem}

\begin{pf} By Lemma 21.2 in van der Waart \cite{VaW},
$F_{m_k}^{-1}(x)$ converges to $F_{m}^{-1}(x)$ for all $x$ at
which $F_{m}^{-1}$ is continuous. Since $F_{m}^{-1}$ is
increasing, this fails to happen for at most countably many
points $x\in[0,1]$. The statement follows.
\end{pf}

\begin{pf*}{Proof of Theorem \ref{theo:principalneq}}
Recall that we always take $\pi=Id$. From Proposition~\ref{prop:bodell}
we have
\begin{eqnarray*}
B_n(n\mu t,\lambda/n)= \mathbb{E} \Biggl(
\sum_{i=1}^n p_i e^{-n\mu t}\prod_{j=1,j\not=i}^n
 [\mathbf{1}_{i<j} +(e^{ n\mu p_j
t}-\mathbf{1}_{i<j})e^{-\lambda/n} ]  \Biggr).
\end{eqnarray*}
Since $\sum_{j} p_j =1$, we can rewrite
\begin{eqnarray*}
&&B_n(n\mu t,\lambda/n)\\
&&\qquad=\mathbb{E} \Biggl( \sum_{i=1}^n p_i e^{-n\mu p_i
t}\prod_{j=1,j\not=i}^n  [1-(1-e^{-\lambda/n})(1- \mathbf{1}_{i<j}
e^{- n\mu p_j t})  ]  \Biggr).
\end{eqnarray*}

Let us define
\[
\tilde{B}_n:=\mathbb{E} \Biggl( \frac{1}{n \mu} \sum_{i=1}^n n \mu p_i
e^{-n \mu p_i t} \exp \Biggl\{-\lambda/n\sum_{j=1,j\not=i
}^{n}(1-e^{-n \mu p_j t}\mathbf{1}_{i<j} ) \Biggr\} \Biggr).
\]

It is elementary to check that $| B_n(n\mu t,\lambda/n) -
\tilde{B}_n|\leq\frac{C}{n}$, so we shall study the term
$\tilde{B}_n$. We have that
\begin{eqnarray*}
\tilde{B}_n&=& \mathbb{E} \Biggl( \exp \Biggl\{-\lambda/n\sum_{j=1}^n
1-e^{-n \mu p_j t} \Biggr\}\\
&&\hspace*{13pt}\times
\frac{1}{n \mu} \sum_{i=1}^n n \mu p_i e^{-n \mu p_i t}
\exp \Biggl\{-\lambda/n\sum_{j=1}^{i} e^{-n \mu p_j
t} \Biggr\}e^{\lambda/n} \Biggr).
\end{eqnarray*}
Therefore, thanks to the bound $xe^{-xt}\leq\frac{1}{t}$ we have
\[
|e^{-\lambda/n}\tilde{B}_n-L(\mu,t,\lambda)|\leq
\frac{1}{\mu\lambda}\mathbb{E} |\Psi(\nu^n) | +
 | \mathbb{E}(\hat{\mathcaligr L}_n(\mu,t,\lambda)) - \hat
{\mathcaligr
L}(\mu,t,\lambda) |,
\]
with
\begin{eqnarray*}
\Psi(m)&:=& \exp \biggl\{-\lambda\int_{\mathbb{R}_+} 1-e^{-xt}
m(dx) \biggr\}-\exp\bigl\{-\lambda\bigl(1-\phi(t)\bigr)\bigr\},
\\
\hat{\mathcaligr L}_n(\mu,t,\lambda)
&:=& \frac{1}{n \mu} \sum
_{i=1}^n n
\mu p_i
e^{-n \mu p_i t} \exp \Biggl\{-\lambda/n\sum_{j=1}^i e^{-n \mu p_j
t} \Biggr\}
\end{eqnarray*}
and $\hat{\mathcaligr L}(\mu,t,\lambda)$ defined as follows:
\begin{eqnarray*}
\hat{\mathcaligr L}(\mu,t,\lambda)&=&\frac{-\phi'(t)}{\mu} \int_0^1
e^{-\lambda
\phi(t)x}\, dx \qquad \mbox{if LLN-$P${-ex} holds,}
\\
\hat{\mathcaligr L}(\mu,t,\lambda)&=&\frac{1}{\mu} \int_0^{\infty} x
e^{-x t}
\exp\biggl\{ -\lambda\int_x^{\infty} e^{- yt} P(dy)\biggr\} P(dx)\qquad
\mbox{if  LLN-$P^-$ holds}
\end{eqnarray*}
or
\[
\hat{\mathcaligr L}(\mu,t,\lambda)=\frac{1}{\mu} \int_0^{\infty} x
e^{-x t}
\exp\biggl\{ -\lambda\int_0^x e^{- yt} P(dy)\biggr\} P(dx)\qquad \mbox{if
 LLN-$P^+$ holds}.
\]
Since $\Psi$ is continuous and bounded in ${\mathcaligr P}(\mathbb
{R}_+)$ and
$\Psi(P)=0$, we get by LLN-$P$ and Theorem~\ref{teoLLN} that
$\mathbb{E} |\Psi(\nu^n) |\to0$ when $n\to\infty$.
Thus, we
just have to prove that
\[
\mathbb{E}(\hat{\mathcaligr L}_n(\mu,t,\lambda))\longrightarrow
\hat
{\mathcaligr L}(\mu,t,\lambda).
\]

\textit{The exchangeable case.} Notice that under LLN-$P${-ex},
\begin{eqnarray*}
&&\mathbb{E}(\hat{\mathcaligr L}_n(\mu,t,\lambda))\\
 &&\qquad= \mathbb{E}
\Biggl( \frac{1}{n \mu}
\sum_{i=1}^n \frac{1}{n!} \sum_{\sigma\in\Pi} n \mu
p_{\sigma(i)} e^{-n \mu p_{\sigma(i)}
t}\\
&&\qquad\quad\hspace*{79pt}\times
\exp \Biggl\{-\lambda/n\sum_{j=1}^{i} e^{-n \mu
p_{\sigma(j)} t} \Biggr\} \Biggr) \\
&&\qquad=
\mathbb{E} \Biggl( \frac{1}{n \mu} \sum_{i=1}^n \sum_{k=1}^n
\frac{1}{n!} \sum_{\sigma\in\Pi, \sigma(i)=k} n \mu
p_{\sigma(i)} e^{-n \mu p_{\sigma(i)}
t}\\
&&\qquad\quad\hspace*{123pt}\times\exp \Biggl\{-\lambda/n\sum_{j=1}^{i} e^{-n \mu
p_{\sigma(j)} t} \Biggr\} \Biggr) \\
&&\qquad=
 \mathbb{E} \Biggl( \frac{1}{n \mu} \sum_{k=1}^n n \mu p_k e^{-n
\mu p_k t}\\
&&\qquad\quad\hspace*{47pt}\times \frac{1}{n} \sum_{i=1}^n \frac{1}{(n-1)!}\sum_{\sigma
\in\Pi, \sigma(i)=k} \exp \Biggl\{-\lambda/n\sum_{j=1}^{i} e^{-n
\mu p_{\sigma(j)} t} \Biggr\} \Biggr) .
\end{eqnarray*}
Since by LLN-$P$ and Theorem \ref{teoLLN},
\begin{eqnarray*}
\mathbb{E} \Biggl( \frac{1}{n \mu} \sum_{k=1}^n n \mu p_k e^{-n \mu p_k
t}  \Biggr)
\int_0^1 e^{-\lambda\phi(t)x}\,dx
\longrightarrow
\frac{-\phi'(t)}{\mu}
\int_0^1 e^{-\lambda\phi(t)x}\,dx
\end{eqnarray*}
when $n\to\infty$, it is enough to show that
$\delta_n(\mu,t,\lambda)$ goes to $0$ when $n\to\infty$, where
\begin{eqnarray*}
\delta_n(\mu,t,\lambda)&:= & \mathbb{E}(\hat{\mathcaligr L}_n(\mu
,t,\lambda
))- \mathbb{E} \Biggl(
\frac{1}{n \mu} \sum_{k=1}^n n \mu p_k e^{-n \mu p_k t}  \Biggr)
\int_0^1 e^{-\lambda\phi(t)x}\,dx \\
&&\hspace*{-10.8pt}= \mathbb{E}(\hat{\mathcaligr L}_n(\mu,t,\lambda))- \frac{1}{\mu
} \mathbb{E}
\Biggl( \int_{\mathbb{R}_+} x e^{-xt} \nu^{(n)}(dx)  \Biggr)
\int_0^1 e^{-\lambda\phi(t)x}\,dx.
\end{eqnarray*}

Let us write for $i=1,\ldots,n-1$, and a permutation $\sigma$ of
$\{1,\ldots,n\}$,
\[
\alpha_t^{\sigma}(i,n):=\sum_{j=1}^{i} e^{-n\mu p_{\sigma(j)}
t}\quad\mbox{and}\quad\alpha_t(n,n):=\sum_{j=1}^{n} e^{-n\mu p_{j}t}.
\]
Define furthermore
\begin{eqnarray*}
I_n^k&=& \frac{1}{n} \sum_{i=1}^n
 \biggl(\exp \biggl\{-\frac{\lambda}{n}
\alpha_t(n,n)\frac{i}{n} \biggr\}\\
&&\hspace*{29pt}\times
\frac{1}{(n-1)!} \sum
_{\sigma\in
\Pi, \sigma(i)=k}  \biggl[ \exp \biggl\{-\lambda/n  \biggl[
\alpha_t^{\sigma}
(i,n)
-\frac{i}{n} \alpha_t(n,n)  \biggr]  \biggr\} -1 \biggr]  \biggr),\\
\mathit{II}_n
&=&
\frac{1}{n} \sum_{i=1}^n \exp \biggl\{-{\frac{\lambda}{n}}
\alpha_t(n,n)\frac{i}{n} \biggr\} - \exp \biggl\{-\lambda\phi(t)
\frac{i}{n} \biggr\}
\end{eqnarray*}
and
\[
\mathit{III}_n=\frac{1}{n} \sum_{i=1}^n \exp \biggl\{-\lambda\phi(t)
\frac{i}{n} \biggr\} - \int_0^1 e^{-\lambda\phi(t)x}\,dx.
\]

Then, we have
\begin{eqnarray*}
|\delta_n(\mu,t,\lambda)|&\leq&
 \Bigg| \mathbb{E} \Biggl( \frac{1}{n \mu} \sum_{k=1}^n n \mu p_k
e^{-n \mu p_k
t}  ( I_n^k+\mathit{II}_n+\mathit{III}_n  )  \Biggr) \Bigg | \\
&\leq& \frac{1}{t \mu}  \Biggl[\frac{1}{n } \sum_{k=1}^n \mathbb{E}|
I_n^k |+ \mathbb{E}|\mathit{II}_n|+ |\mathit{III}_n|  \Biggr]
\end{eqnarray*}
thanks to the bound $x e^{-xt}\leq\frac{1}{t}$.
Term $\mathit{III}_n$ clearly goes to $0$ when
$n\to\infty$. On the other hand, we have
\begin{eqnarray*}
\mathbb{E}|\mathit{II}_n|
\leq\frac{\lambda}{n} \sum_{i=1}^n
\frac{i}{n}\mathbb{E} \bigg|\frac{1}{n}\alpha_t(n,n)-\phi(t)
\bigg|
\leq
\frac{\lambda}{2} \mathbb{E} \bigg|\int_{\mathbb{R}_+} e^{- x t}
\nu^n(dx)-\phi(t)\bigg |.
\end{eqnarray*}
The mapping $\nu\mapsto|\int_{\mathbb{R}_+} e^{- x t} \nu(dx)
-\int_{\mathbb{R}_+} e^{- x t} P(dx)|$ being continuous and bounded on
${\mathcaligr P}(\mathbb{R}_+)$, the latter term goes to $0$ by LLN-$P$
and Theorem~\ref{teoLLN}.

Now, by exchangeability $\mathbb{E}| I_n^k |$ does not depend on $k$, and
moreover, setting $\alpha_t(i,n):=\sum_{j=1}^{i} e^{-n\mu p_{j}
t}$, we have
\begin{eqnarray*}
\mathbb{E}|I_n^k|
&\leq&
 \frac{1}{n} \sum_{i=1}^n \mathbb{E}
\bigg|\exp
 \biggl\{-\lambda/n  \biggl[ \alpha_t(i,n)
-\frac{i}{n} \alpha_t(n,n)  \biggr]  \biggr\}-1  \bigg| \\
&\leq&  \frac{1}{n} \sum_{i=1}^n \mathbb{E} \Bigg|\lambda/n \Biggl [
\sum_{j=1}^i e^{-n \mu p_j t} -\frac{i}{n} \sum_{k=1}^n e^{-n \mu
p_k t} \Biggr]
 \Bigg| \\
&= & \frac{1}{n} \sum_{i=1}^n \mathbb{E} \Bigg|\lambda/n
\sum_{j=1}^i  \biggl( e^{-n \mu p_j t}- \int e^{- x t}
\nu^{(n)}(dx)  \biggr) \Bigg |,
\end{eqnarray*}
and so
\[
\mathbb{E}|I_n^k|
\leq \frac{1}{n} \sum_{i=1}^n \mathbb{E}
\Bigg|\lambda/n
\sum_{j=1}^i  \bigl( e^{-n \mu p_j t}- \phi(t)  \bigr)  \Bigg|
+
\frac{\lambda}{2} \mathbb{E} \bigg|\phi(t)- \int e^{- x t} \nu
^{(n)}(dx) \bigg|.
\]
Thus, we just have to check that $\mathit{IV}_n:=\frac{\lambda}{n^2}
\sum_{i=1}^n \mathbb{E} |
\sum_{j=1}^i  ( e^{-n \mu p_j t}- \phi(t)  )  |$
goes to $0$. Indeed, we have
\begin{eqnarray*}
\mathit{IV}_n&\leq& \frac{\lambda}{n^2} \sum_{i=1}^n \Biggl [\mathbb{E} \Biggl(
\sum_{j=1}^i  \bigl( e^{-n \mu p_j t}- \phi(t) \bigr)  \Biggr)^2
 \Biggr]^{1/2} \\
&= & \frac{\lambda}{n^2} \sum_{i=1}^n \Biggl [
\sum_{j=1}^i \mathbb{E} \bigl( e^{-n \mu p_j t}- \phi(t) \bigr)^2
\\
&&\hspace*{34pt}{}+ \sum_{k=1}^i \sum_{l=1,l\not=k }^i \mathbb{E}\bigl( e^{-n
\mu p_l t}- \phi(t)\bigr)\bigl( e^{-n \mu p_k t}- \phi(t)\bigr)
 \Biggr]^{1/2} \\
&\leq&  \frac{2\lambda}{\sqrt{n}} + \frac{\lambda}{n^2}
\sum_{i=1}^n  \bigl[ i(i-1) \big| \mathbb{E}\bigl( e^{-n \mu p_1 t}-
\phi(t)\bigr)\bigl( e^{-n \mu p_2 t}- \phi(t)\bigr) \big|
 \bigr]^{1/2}. \end{eqnarray*}
Therefore,
\begin{eqnarray*}
\mathit{IV}_n \leq&  \dfrac{2\lambda}{\sqrt{n}} + \lambda
 \big| \mathbb{E}\bigl( e^{-n \mu p_1 t}- \phi(t)\bigr)\bigl( e^{-n \mu p_2 t}-
\phi(t)\bigr)
 \big|^{1/2}.
\end{eqnarray*}
By LLN-$P$ and Proposition \ref{propchaos}(i) with
$k=2$, we conclude that the latter term goes to $0$. This finishes
the proof in the exchangeable case.

\textit{The monotone cases.} We consider the case when  LLN-$P^+$ holds, the decreasing case being similar. Notice that
if $F_n^{-1}(x):=\inf\{t\geq0\dvtx F_n(t)\geq x\}$ is the
generalized inverse of $F_n(x)=\nu^{(n)}([0,x])$, we have that
\[
\hat{\mathcaligr L}_n(\mu,t,\lambda):= \frac{1}{ \mu} \int_0^1
F_n^{-1}(x)e^{-F_n^{-1}(x)t} \exp \biggl\{-\lambda\int_0^{i_n(x)}
e^{-F_n^{-1}(y) t}\, dy  \biggr\} \,dx ,
\]
where $ i_n(x)=\frac{ \lceil nx \rceil}{n}$ and $\lceil\cdot
\rceil$ is the ceiling function. On the other hand, under the law
$dx$ the generalized inverse $F^{-1}\dvtx[0,1]\to\mathbb{R}_+$ of $F$ is a
random variable of law~$P$. We thus have
\[
\hat{\mathcaligr L}(\mu,t,\lambda):= \frac{1}{ \mu} \int_0^1
F^{-1}(x)e^{-F^{-1}(x)t} \exp \biggl\{-\lambda\int_0^{x}
e^{-F^{-1}(y) t}\, dy  \biggr\} \,dx.
\]

Thanks to this and the bound $x e^{-xt}\leq\frac{1}{t}$, we get
that
\begin{eqnarray*}
&&|\mathbb{E}(\hat{\mathcaligr L}_n(\mu,t,\lambda))-\hat{\mathcaligr
L}(\mu
,t,\lambda)|\\
&&\qquad\leq
 \frac{\lambda}{ t \mu} \mathbb{E} \biggl(\int_0^1 \int_0^{x}
\big|e^{-F_n^{-1}(y)t}-e^{-F^{-1}(y)t}\big|\, dy\, dx \biggr) \\
&&\quad\qquad{} + \frac{\lambda}{ t \mu} \int_0^1 \int_x^{i_n(x)}
e^{-F^{-1}(y)t}\, dy\, dx \\
&&\quad\qquad{} +
 \frac{1}{ \mu} \mathbb{E} \biggl( \int_0^1
 \big|F_n^{-1}(x)e^{-F_n^{-1}(x)t}-F^{-1}(x)e^{-F^{-1}(x)t}\big |
\,dx \biggr)\\
&&\qquad \leq
\frac{\lambda}{ t \mu} \mathbb{E} \biggl(\int_0^1
\big|e^{-F_n^{-1}(y)t}-e^{-F^{-1}(y)t}\big|\, dy  \biggr) + \frac{\lambda}{ n
t \mu
} \\
&&\quad\qquad{} + \frac{1}{ \mu} \mathbb{E} \biggl( \int_0^1
 \big|F_n^{-1}(x)e^{-F_n^{-1}(x)t}-F^{-1}(x)e^{-F^{-1}(x)t}\big |
\,dx \biggr).
\end{eqnarray*}
Therefore, and thanks also to LLN-$P$ and Theorem~\ref{teoLLN},
it is enough to prove that the bounded functionals
on ${\mathcaligr P}(\mathbb{R}_+)$
\[
\nu\mapsto\int_0^1 \big|e^{-F_{\nu}^{-1}(y)t}-e^{-F^{-1}(y)t}\big|\, dy
\]
and
\[
\nu\mapsto\int_0^1
 \big|F_{\nu}^{-1}(x)e^{-F_{\nu
}^{-1}(x)t}-F^{-1}(x)e^{-F^{-1}(x)t} \big|\,dx
\]
are continuous, since they both vanish at $\nu=P$. This follows
by dominated convergence and Lemma~\ref{lemma:continvdistrfunc}.
The proof of the theorem is finished.
\end{pf*}

\section{Concluding remarks}\label{sec5}

The limiting stationary regime of the MtF search-cost as the
number of objects tend to infinity has been considered by several
authors. One of the motivations is to compare efficiency among
different popularity distributions when equilibrium is reached.
Nonetheless, the rate at which equilibria are reached should also
account for efficiency considerations. This was one of the
motivations of the present article.

We have developed a general framework for studying the limiting
dynamical behavior of the MtF search-cost when requests rates are
sampled from empirical probability measures that asymptotically
approach a specified law $P$. In this law of large numbers
asymptotic regime, popularities of objects are comparable, in the
sense that their asymptotic average is finite and nonnull. By
this reason, although the transient behavior depends on the
initial ordering, is not considerably sensitive to it. This can be
seen in the fact that a common convergence rate to equilibrium
$O(\int x e^{-tx}P(dx))$ was obtained in the three representative
situations considered.\looseness=1

Our techniques also ensure the asymptotic stability under
perturbations of request rates that preserve $P$. Namely, the
limiting expectations of functionals of $(p_i^{(n)})_{i=1}^n$,
which are symmetric functions of $\nu^{(n)}$, depend only on $P$.
This was the case of the equilibrium part of the transient
search-cost. The out-of-equilibrium transient search-cost involved
in turn nonsymmetric functionals of $(p_i^{(n)})_{i=1}^n$, which
yielded different limits according to the different
``enumerations'' of objects. This is the explanation for the
coincidences and discrepancies pointed out in examples~(5) and (6)
of Section~\ref{sec3}. Nevertheless, under the assumption of
exchangeability one still might replace nonsymmetric functionals
by symmetrized versions of them (as in the proof of Theorem~\ref{theo:principalneq})
and obtain ``symmetric'' limits.

It is in principle possible to use our techniques in the
asymptotic analysis of other sorting algorithms, at least in those
cases where the corresponding relevant variables depend on the
empirical measures of the popularities or of the request rates.
However it is
not obvious to identify which functionals of the empirical
measures are involved.

On the other hand, the law of large numbers asymptotic behavior we
have described corresponds to a very particular scaling limit, in
the sense explained before. Therefore, it a priori excludes
deterministic cases of interest such as the Zipf laws
$w_i=i^{\alpha}$ with $\alpha\leq-1$ or scaling approximations
of the Poisson--Dirichlet distribution (see, e.g., \cite{Kingman}, Chapter~9,
 and Joyce and Tavar\'{e} \cite{JT} for the limits of symmetric
linear functionals of these random partitions). Neither the fluid
limit approximation of the search-cost studied by Jelenkovi\'{c}
\cite{Jelen} is covered in its whole generality by our approach.
Indeed, if $Q$ has unbounded support, one might try to approximate
$Q$ by compactly supported laws as in point (d) of Section~\ref{sec3}.
But if one chooses therein $c=c_n$ diverging with $n$, the
empirical means vanish as $n$ goes to infinity.

In these examples, the transient dynamics may be of particular
interest, since they exhibit coexistence of microscopic and macroscopic
popularities which is likely to affect the convergence to
equilibrium. A similar question could be of interest in the
context of the PAC algorithm introduced in \cite{JR}. The
splitting of the transient search-cost we introduced here could be
useful in those cases. A combination of ideas in
\cite{BarreraHuilletParoissinMtF2} and results in \cite{JT}
could help to extend part of our arguments to Poisson--Kingman type
asymptotics, although additional difficulties arise. The
computation of the search-cost law involved highly nonlinear
functionals of the empirical measures, which cannot be deduced
from the asymptotics results obtained in \cite{JT}. These and
related questions are addressed in progressing works by the
authors.

\section*{Acknowledgments}

The authors thank the anonymous referee
for drawing our attention to the papers \cite{Jelen} and
\cite{JR}, so as for several suggestions that allowed us to
improve the presentation of this work.

%

\printaddresses

\mbox{}

\end{document}